\documentstyle
[pb-diagram,lamsarrow,pb-lams,epic,eepic,amssymb,amscd,11pt, fullpage]
{amsart}

\newtheorem{theorem}{Theorem}[section]
\newtheorem{Theorem}[theorem]{Theorem}
\newtheorem{definition}[theorem]{Definition}

\newtheorem{lemma}[theorem]{Lemma}
\newtheorem{Lemma}[theorem]{Lemma}
\newtheorem{corollary}[theorem]{Corollary}

\newtheorem{proposition}[theorem]{Proposition}
\newtheorem{example}[theorem]{Example}

\newtheorem{remark}[theorem]{Remark}
\newtheorem{Remark}[theorem]{Remark}

\newcommand{\pr} {\smallskip\noindent{\bf Proof\,\,}}
\newenvironment{Proof}	{\pr}{\hspace*{\fill}\qed\\}
\newenvironment{proof}	{\pr}{\hspace*{\fill}\qed\\}

\newcommand\bbR{{\Bbb R}}
\newcommand\R{{\Bbb R}}
\newcommand\bbC{{\Bbb C}}
\newcommand\C{{\Bbb C}}
\newcommand\bbZ{{\Bbb Z}}
\newcommand\Z{{\Bbb Z}}
\newcommand\bbT{{\Bbb T}}

\newcommand\cm{{\cal M}}

\newcommand\cC{{\cal C}}
\newcommand\cp{{\cal P}}

\newcommand\fg {{\frak g}}
\newcommand\fh {{\frak h}}
\newcommand\fk {{\frak k}}
\newcommand\ft {{\frak t}}

\newcommand\cS {{\cal S}}
\newcommand\U {{\cal U}}
\newcommand\tcU {{\tilde{\cal U}}}
\newcommand\tcV {{\tilde{\cal V}}}

\newcommand\intD{\mathaccent23\Delta}
\newcommand\intF{\mathaccent23{F}}

\def	\tU	{{ \tilde{U}  }}
\def	\tX	{{ \tilde{X}  }}
\def	\tG	{{ \tilde{G}  }}
\def	\tV	{{ \tilde{V}  }}
\def	\tT	{{ \tilde{T}  }}
\newcommand{\im}	{\operatorname{im}}

\newcommand{\Gl}	{\operatorname{GL}}
\newcommand{\SP}	{\operatorname{Sp}}
\newcommand{\Su}	{\operatorname{SU}}
\newcommand{\lra}	{\longrightarrow}

\begin{document}

\title[Torus actions on symplectic orbifolds \quad July 23, 1995
]{Hamiltonian torus actions on symplectic orbifolds and toric
varieties}

\author[E. Lerman and S. Tolman]{Eugene Lerman and Susan
Tolman \quad July 23, 1995}

\address{Department of Mathematics, MIT, Cambridge, MA 02139}

\curraddr{Department of Mathematics, Univ.\ of Illinois at
Urbana-Champaign, 1409 West Green St., Urbana, IL 61801}
\email{eugene@@math.mit.edu, lerman@@math.uiuc.edu}
\address{Department of Mathematics, MIT, Cambridge, MA 02139}
\email{tolman@@math.mit.edu}
\thanks{Both authors were partially supported by NSF
postdoctoral fellowships.}



\maketitle
\begin{abstract}
In the first part of the paper, we build a foundation for further work
on Hamiltonian actions on symplectic orbifolds.  Most importantly we
prove the orbifold versions of the abelian connectedness and convexity
theorems.

In the second half, we prove that compact symplectic orbifolds with
completely integrable torus actions are classified by convex simple
rational polytopes with a positive integer attached to each facet and
that all such orbifolds are algebraic toric varieties.
\end{abstract}

\tableofcontents
\section{Introduction}\label{section_intro}

\noindent
In this paper we study Hamiltonian torus actions on symplectic orbifolds,
with an emphasis on completely integrable actions.

The category of Hamiltonian group actions on symplectic manifolds is
not closed under symplectic reduction; generically, the reduced space
is an orbifold.  In contrast the symplectic reduction of an orbifold
is generically an orbifold.  Symplectic reduction is a powerful
technique which has been used successfully in such diverse areas as
Hamiltonian systems and representation theory.

Therefore we need to understand symplectic orbifolds even if we only
want to understand symplectic manifolds.
For example, a proof of the Guillemin-Sternberg conjecture that
quantization commutes with reduction naturally encounters orbifolds
\cite{DGMW}, \cite{Meinrenken}.  In the same spirit, the proof of the
non-abelian convexity theorem for manifolds can be reduced to an
abelian convexity theorem for orbifolds \cite{LMTW}.  Orbifolds also
arise in the study of resonances in Hamiltonian systems.

Their are two  main purpose of this paper.

In the first half of the paper, we build a foundation for further work on
Hamiltonian actions on symplectic orbifolds.  For example, we classify
the neighborhoods of isotropic orbits, and we discuss the extension of
Bott-Morse theory to orbifolds.  Most importantly, we prove the
following analogues of the abelian connectedness and convexity
theorems (see \cite{A} and \cite{GS}).\\[5pt]

\noindent {\bf Theorem } {\em
Let a torus $T$ act on a compact symplectic orbifold $(M,\omega)$,
with a moment map $\phi:M \to \ft^*$.
For every $a \in \ft^*$, the fiber $\phi^{-1}(a)$ is connected.
}\\

\noindent {\bf Theorem } {\em
Let a torus $T$ act on a compact symplectic orbifold $(M,\omega)$,
with a moment map $\phi:M \to \ft^*$.
The image of the moment map $\phi(M) \subset \ft^*$ is a rational convex
polytope.  In particular it is the convex hull of the image of the
points in $M$ fixed by $T$,
$$
	\phi (M) = \text{convex hull } (\phi (M^T)).
$$
}

In the second half of the paper, we consider the special case of
completely integrable torus actions, and prove the following theorem:
\\

\noindent {\bf Theorem } {\em
Compact symplectic toric orbifolds
are classified by convex rational simple polytopes with a positive
integer attached to each facet.
}\\

This theorem generalizes a theorem of Delzant \cite{Del} to the case
of orbifolds.  He proved that symplectic toric {\em manifolds} are
classified by the image of their moment maps, that is, by a certain
class of rational polytopes.  It is easy to see that additional
information is necessary for orbifolds:
\begin{example}{\em
Given positive integers $n$ and $m$, there exists a symplectic
orbifold $M$ which is topologically a two sphere, but which look
locally like $\C/(\Z/n\Z)$ and $\C/(\Z/m\Z)$ near its north a south
pole, respectively.  The circle action which rotates $M$ around its
north-south axis is Hamiltonian.  Although the image of the moment map
is a line interval for all $m$ and $n$, these orbifolds are not
isomorphic.  }\end{example}

To state the above theorem precisely. we must  define a few terms.

\begin{definition}
\label{def zero}
A {\em symplectic toric orbifold}
is a quadruple $(M, \omega, T, \phi)$ where $\omega$ is a symplectic
form on a  connected orbifold $M$ and $\phi :M \to \ft^*$ is a
moment map for an effective Hamiltonian action of a torus $T$ on $M$
such that $\dim T = \frac{1}{2} \dim M$.\footnote
{This is not the definition used in algebraic geometry.
We will show that every
compact symplectic toric orbifold can be given the structure of an algebraic
toric variety.}
Two symplectic toric orbifolds are {\em isomorphic} if they are
equivariantly symplectomorphic (implicitly the torus is fixed in this
definition).
\end{definition}

\begin{definition}
\label{def one}
Let $x$ be a point in an orbifold $M$, and let $(\tU,\Gamma,
\phi)$ be a uniformizing chart for neighborhood $U$ of $x$ (see
\cite{Sa} or section~\ref{section_action}), then the {\em (orbifold)
structure} group of $x$ is the isotropy group of $\tilde{x} \in \tU$,
where $\phi(\tilde{x}) = x$.  This group is well defined as an
abstract group.
\end{definition}

\begin{definition}
\label{def two}
Let $\ft$ be a vector space with a lattice $\ell$; let $\ft^*$ denote
the dual vector space.
A convex polytope $\Delta \subset \ft^*$ is {\em rational} if
$$
\Delta = \cap_{i=1}^N \{\alpha \in \ft^* \mid \left< \alpha, y_i
\right> \geq \eta_i \}
$$
for some $y_i \in \ell$ and $\eta_i \in \R$. A {\em facet} is a face
of codimension one.  An $n$ dimensional polytope is {\em simple} if
exactly $n$ facets meet at every vertex.  For this paper, we shall
adopt the convenient but non standard abbreviation that a {\em labeled
polytope} is a convex rational simple polytope plus a positive integer
attached to each facet.  Two labeled polytopes are {\em isomorphic} if
one can be mapped to the other by a translation and the corresponding
facets have the same integer labels.
\end{definition}

\begin{theorem}
\label{classification}
For a compact symplectic toric orbifold $(M,\omega,T, \phi )$ the
image of the moment map $\phi (M)$ is a rational simple polytope.  For
every facet $F$ of $\phi (M)$ there exists a positive integer $n_F$
such that for every $x$ in the preimage of the interior of the facet
the structure group of $x$ is $\Z/n_F\Z$.

Two compact symplectic toric orbifolds are isomorphic if and only if
their associated labeled polytopes are isomorphic.

Every labeled polytope occurs as the image of a compact symplectic
toric orbifold.
\end{theorem}

Finally, we consider K\"ahler structures on symplectic toric
orbifolds.

\begin{definition}
\label{def four}
Let $\Delta \subset \ft^*$ be a convex polytope  such that $\dim(\Delta)
= \dim(\ft)$.
Given a face $F$ of $\Delta$, the {\em dual cone} to
$F$ is the set
$$
C_F = \{ \alpha \in \ft \mid \left< \alpha, \beta -\beta'\right> \leq
0 \text{ for all } \beta \in F \text{ and } \beta' \in \Delta\}.
$$
The {\em fan} of $\Delta$  is the set of dual cones to the faces of
$\Delta$.
\end{definition}

\begin{theorem}
Every compact symplectic toric orbifold admits an equivariant
complex structure which is compatible with its symplectic form.
Two symplectic toric orbifolds with compatible  complex structures
are equivariantly biholomorphic exactly if their image polytopes
have the same fans.
\end{theorem}

\noindent {\sc Acknowledgments}\\
\vspace{-12pt}

At the conference {\em Applications of Symplectic Geometry} at the
Newton Institute, 10/31/94 - 11/11/94, we learned that R. de Souza and
E. Prato have been working independently on the same problem.

It is a pleasure to thank Chris Woodward for many useful discussions.
In particular, section~\ref{local to global} is a joint work with
Chris Woodward.  We thank Sheldon Chang for a number of helpful
suggestions.  \\

\part{Hamiltonian torus actions on symplectic orbifolds}

\section{Group actions on orbifolds}\label{section_action}

In this section, we recall definitions related to actions of groups on
orbifolds and describe some properties of actions of compact groups.
The main result is the slice theorem.  The presentation is largely
self-contained and borrows heavily from a paper of Haefliger and Salem
\cite{HS}.

We begin by defining orbifolds and related differential geometric notions.

An {\em orbifold} $M$ is a topological space $|M|$, plus an {\em
atlas} of {\em uniformizing charts} $(\tU,\Gamma,\varphi)$, where
$\tU$ is open subset of $\R^n$, $\Gamma$ is a finite group which acts
linearly on $\tU$ and fixes a set of codimension at least two, and
$\varphi: \tU \to |M|$ induces a homeomorphism from $\tU/\Gamma$ to $U
\subset |M|$.  Just as for manifolds, these charts must cover $|M|$;
they are subject to certain compatibility conditions; and there is a
notion of when two atlases of charts are equivalent.  For more
details, see Satake \cite{Sa}.  Given $x$, we may choose a
uniformizing chart $(\tU,\Gamma,\varphi)$ such that $\varphi^{-1}(x)$
is a single point which is fixed by $\Gamma$.  In this case $\Gamma$
is the orbifold structure group of $x$, cf.\
Section~\ref{section_intro}.

Given a point $x$ in an orbifold $M$ and a uniformizing chart
$(\tU,\Gamma,\varphi)$ with $\varphi^{-1}(x)$ a single point, define
the {\em uniformized tangent space at $x$} to be the tangent space to
$\varphi^{-1}(x)$ in $\tU$ and denote it by $\tT_xM$.  The
quotient $\tT_x M/\Gamma$ is $T_x M$, the fiber of the tangent bundle
of $M$ at $x$.

A vector field $\xi$ on $M$ is a $\Gamma$ invariant vector field
$\tilde{\xi}$ on each uniformizing chart $(\tU,\Gamma,\phi)$; of
course, these must agree on overlaps.  Similar definition apply to
differential forms, metrics, etc.

Let $M$ and $N$ be orbifolds with atlases $\tcU$ and $\tcV$.
A {\em map of orbifolds} $f: M \to N$ is a map $f : \tcU \to \tcV$,
and an equivariant map $\tilde{f}_\tU: (\tU,\Gamma)
\to (\tV,\Upsilon)$ for each $(\tU,\Gamma) \in \tcU$,
where $(\tV,\Upsilon) = f(\tU,\Gamma)$.  These $\tilde{f}_\tU$
are subject to a compatibility condition
which insures, for instance, that $f$ induces a continuous map of
the underlying spaces.  Additionally, there is a notion of
when two such maps are equivalent.  Again, see \cite{Sa} for details.

\begin{definition}
\label{def_action}{\rm
Let $G$ be a Lie group.  A {\em smooth action} $a$ of $G$ on an
orbifold $M$ is a smooth orbifold map $a: G\times M \to M$ satisfying
the usual group laws, that is, for all $g_1, g_2 \in G$ and $x \in M$
$$
a(g_1, a(g_2, x )) = a(g_1 g_2, x ) \quad \hbox{and} \quad
a(1_G, x) = x,
$$
where $1_G$ is the identity element of $G$, and  ``=" means
``are equivalent as maps of orbifolds."
}
\end{definition}
Definition~\ref{def_action} implies that
the action $a$ induces a continuous action $|a|$ of $G$ on the
underlying topological space $|M|$.
In particular, for every  $g_0 \in  G$
and $x_0 \in M$ there are neighborhoods $W$ of $g_0$ in $G$, $U$ of $x_0$
in $M$, and $U'$ of $a(g_0,x_0) $ in $M$, charts
$(\tilde{U}, \Gamma, \varphi) $ and $(\tilde{U}', \Gamma ', \varphi') $ and
a smooth map $\tilde{a}: W\times \tilde{U} \to \tilde{U}'$ such that
$\varphi' (\tilde{a}(g, \tilde{x})) =|a|(g, \varphi (x))$ for all $(g,
\tilde{x}) \in W\times \tilde{U}$.  Note that $\tilde{a}$ is not
unique, it is defined up to composition with elements of the orbifold
structure groups $\Gamma $ of $x_0$ and $\Gamma' $ of $g_0\cdot x_0$.

If $g_0 = 1_G$, the identity of $G$, then we may assume $\tU \subset
\tU'$, and we can choose $\tilde{a}$ such that $\tilde{a}(1_G, x) =
x$. Then $\tilde{a}$ induces a local action of $G$ on $\tU$.

An action of a Lie group $G$ on an orbifold $M$ induces an
infinitesimal action of the Lie algebra $\fg$ of $G$ on $M$.  Denote
for a vector $\xi \in \fg$ the corresponding induced vector field by
$\xi_M$.  In particular, for any chart $(\tU,\Gamma,\varphi)$ there
exists a $\Gamma $ invariant vector field $\xi _{\tU}$ on $\tU$ and
such vector fields satisfy compatibility conditions.

If a point $x$ is {\em fixed} by the action of $G$ and $G$ is
compact, then the local action $G$ on the uniformizing chart $\tU$
generates an action of $\tG$ on a subset $\tV \subset \tU$, where
$\tilde{G}$ is a cover of the identity component of $G$.  Note that
the actions of $\tilde{G}$ and $\Gamma $ on $\tV$ {\em commute};
otherwise a group action would not induce the corresponding
infinitesimal actions of the Lie algebra.

More generally one can show that for a fixed point $x$ with structure
group $\Gamma$ there exists a uniformizing chart
$(\tU,\Gamma,\varphi)$ for a neighborhood $U$ of $x$, an exact
sequence of groups
$$
	1\to \Gamma \to \hat{G} \stackrel{\pi}{\to} G \to 1,
$$
and an action of $\hat{G}$ on $\tU$ such that the following diagram
commutes:
$$
\begin{diagram}
\node{\hat{G} \times \tU} \arrow{e} \arrow{s}\node{\tU}\arrow{s}\\
\node{G \times U }	\arrow{e}	\node{U}
\end{diagram}.
$$
The extension $\hat{G}$ of $G$ depends on $x$ and, in particular, is
not globally defined.  The homomorphism $\pi :\hat{G} \to G$ induces
an isomorphism of the Lie algebras $\varpi: \hat{\fg} \to \fg$ and for
any $\xi \in \hat{\fg}$ we have $\xi_{\tU} = (\varpi (\xi))_{\tU}$.

The simplest examples of group actions on orbifolds are linear actions
of groups on vector orbi-spaces.  A {\em vector orbi-space} is
a quotient of the form
$V/\Gamma$ where $V$ is a vector space and $\Gamma$ is a finite subgroup of
$\Gl (V)$.
We define
$$
	\Gl (V/\Gamma) := N(\Gamma)/\Gamma,
$$
where $N(\Gamma)$ is the normalizer of $\Gamma$ in $\Gl (V)$.
The group $\Gl (V/\Gamma)$ does act on the orbifold $V/\Gamma$ in
the sense of Definition~\ref{def_action}.
We define a {\em representation} $\rho :H \to \Gl (V/\Gamma)$ of a group $H$
on the vector orbi-space $V/\Gamma$ to be a group homomorphism $\rho :
H \to N(\Gamma )/\Gamma$.  A representation of $H$ on
$V/\Gamma$ defines an action of $H$ on the orbifold $V/\Gamma$.
For a more detailed discussion of representations, please see
Section~\ref{secsymploc}.

Let $G$ be a compact Lie group acting on an orbifold $M$.  Let $x$ be
a point in $M$, let $\Gamma$ be its orbifold structure group, and let
$G_x$ be its stabilizer.  Because $\Gamma$ commutes with the local
action of $G_x$ on a chart $\tU$, the uniformized tangent space $\tT_x
(G\cdot x) \subset \tT_xM$ is fixed by $\Gamma$.  Thus, there is a
natural representation of the isotropy group $G_x$ of $x$ on the
vector orbi-space $W/\Gamma$, where $W= \tT_x M / \tT_x (G\cdot x)$
(we may also identify $W$ with the orthogonal complement to $\tT_x
(G\cdot x)$ in $\tT_x M$ with respect to some invariant metric).
Because of Proposition~\ref{slice_prop} below we will refer to the
vector orbi-space $W/\Gamma$ as the {\em (differential) slice} for the
action of $G$ at $x$, and to the representation $G_x \to \Gl
(W/\Gamma)$ as the {\em (differential) slice representation}.

\begin{proposition}\label{slice_prop}
{\rm ({\bf Slice theorem})}
Suppose that a compact Lie group $G$ acts on an orbifold $M$ and $G\cdot x$
is an orbit of $G$. A $G$ invariant neighborhood of the orbit is
equivariantly diffeomorphic to a neighborhood of the zero section in
the associated orbi-bundle $G\times _{G_x} W/\Gamma$, where $G_x$ is
the isotropy group of $x$ with respect to the action of $G$, $\Gamma$
is the orbifold structure group of $x$, and $W = \tT_x M/ \tT_x (G\cdot x)$.
\end{proposition}

\begin{proof}
This is strictly analogous to the slice theorem for actions on manifolds,
and follows from the fact that metrics can be averaged over
compact Lie groups.
\end{proof}

\begin{remark}\label{orbiproof}{\em
Throughout this paper, the reader will find many proofs similar to the
one above, which simply claim that the the proof for orbifolds is
strictly ``analogous" to the proof for manifolds.  By this, we mean
that because the usual proof (or the particular proof cited) is
functorial, it will also work for orbifolds.  Because each local
uniformizing chart is itself a manifold, we can apply the manifold
proof to construct the desired object on it.  If the construction is
natural, these local objects will form a global object on the
orbifold.  Sometimes the construction depends on an additional
structure, but is natural once that structure is chosen.  In this
case, we choose that structure on the orbifold, and then apply the
above reasoning.

For instance one step in the proof of the theorem above is to show that
there exists a neighborhood of the zero section of the
normal bundle of a suborbifold $X \subset M$ which is diffeomorphic to
a neighborhood of $X$ in $M$.
First, we choose a metric on $M$.
Let's examine the naturality condition explicitly in this case.
Let $\tU$ and $\tU'$ be manifolds with metrics and sub-manifolds
$\tX$ and $\tX'$ respectively.
Let $N(\tX)$ and $N(\tX')$ denote the normal bundles of $\tX$ and $\tX'$.
Let $\lambda : \tU \to \tU'$ be an open isometric embedding such that
$\lambda(\tX) = \lambda(\tU) \cap \tX'$.
The embedding $\lambda$ induces a map $\lambda_* : N(\tX) \mapsto  N(\tX')$.
Let $\psi : N(\tX) \to \tU$ and $\psi' : N(\tX') \to  \tU'$ be the
diffeomorphisms constructed in the proof of the tubular neighborhood
theorem for manifolds.  The construction of $\psi$ is natural in
the sense that the following diagram commutes:
$$
\begin{diagram}
\node{N(\tX)} \arrow{e,t}{\psi}
\arrow{s,l}{\lambda_*}\node{\tU}\arrow{s,r}{\lambda}\\
\node{N(\tX)}	\arrow{e,t}{\psi'}	\node{\tU'}
\end{diagram}.
$$
Therefore, these  maps $\psi$ form a diffeomorphism of orbifolds from
a neighborhood of the zero section of the normal bundle of the suborbifold
to a neighborhood of the suborbifold.

Most authors prefer natural constructions and make scrupulously
clear which choices are necessary.  Therefore, for the most part,
we have not found it necessary to repeat what has been done
well elsewhere.
}\end{remark}

\begin{remark}{\em
As in the case of manifolds, the compactness of the group $G$ is not
necessary for the slice theorem.  It is enough to require that
the induced action on the underling topological space is proper.}
\end{remark}

For an action of a connected group $G$ on an orbifold $M$, it follows
from the existence of slices that the fixed point set $M^G$ is a
suborbifold.\footnote{ Strictly speaking, $M^G$ might have codimension
1 strata and so would not be an orbifold in the sense of this paper.
However, for symplectic actions on symplectic orbifolds, the
codimension of $M^G$ is at least two.}  \, Therefore, the
decomposition of $M$ into infinitesimal orbit types is a
stratification into suborbifolds.  On the other hand, the fixed point
set $M^G$ need not be a suborbifold if the group $G$ is not connected.

\begin{example}{\em
Let $\Gamma = \Z/(2\Z)$ act on $\C^2$ by sending $(x,y)$ to $(-x,-y)$
Let $G = \Z/(2\Z)$ act on $\C^2/\Gamma$ by sending $[x,y]$ to
$[x,-y]$.  Then the fixed point set $M^G$ is $\{[x,0]\} \cup
\{[0,y]\} = \C/\Gamma \cup \C/\Gamma$.}
\end{example}
Consequently the decomposition of an orbifold according to the orbit
type is not a stratification.  Fortunately the following lemma
holds.

\begin{lemma}\label{lem.princ.orbittype}
If $G$ is a compact Lie group acting on a connected orbifold $M$ then there
exists an open dense subset of $M$ consisting of points with the same orbit
type.
\end{lemma}
\begin{proof}
We first decompose the orbifold into the open dense set of smooth
points $M_{\text{smooth}}$ and the set of singular points. Since we
assume that all the singularities have codimension 2 or greater,
$M_{\text{smooth}}$ is connected. A smooth group action preserves this
decomposition. Since $G$ is compact, the action of $G$ on
$M_{\text{smooth}}$ has a principal orbit type (see for example
Theorem~4.27 in \cite{Kaw}). The set of points of this orbit type is
open and dense in $M_{\text{smooth}}$, hence open and dense in $M$.
\end{proof}

\begin{corollary}\label{cor.loc_free}
If a torus $T$ acts effectively on a connected orbifold $M$ then the action
of $T$ is free on a dense open subset of $M$.
\end{corollary}

\section{Symplectic local normal forms}
\label{secsymploc}

In this section, we write down normal forms for the neighborhoods of
isotropic orbits of a compact Lie group $G$ which acts on a symplectic
manifold $(M,\omega)$ in a Hamiltonian fashion; that is, we classify
such neighborhoods up to $G$ equivariant symplectomorphisms.  We also
point out consequences of these normal forms.

The definitions of symplectic manifolds, symplectic group actions,
moment maps, and Hamiltonian actions carry over verbatim to the
category of orbifolds.  To wit, a {\em symplectic orbifold} is an
orbifold $M$ with a closed non-degenerate $2$-form $\omega$.  A group
$G$ acts {\em symplectically} on $(M,\omega)$ if the action preserves
$\omega$.  A moment map $\phi: M \to \fg^*$ for this action is a $G$
equivariant map such that
$$
	\iota(\xi_M )\, \omega = d\langle \xi, \phi\rangle,
	\text{ for all } \xi \in \fg.
$$
If there is a moment map, we say that $G$ acts on $(M,\omega)$ in a
{\em Hamiltonian fashion}.

The simplest symplectic orbifold is a {\em symplectic vector
orbi-space} $V/\Gamma$, where $V$ is a symplectic vector space and
$\Gamma$ is a finite subgroup of the symplectic group
$\SP(V)$.\footnote {The following discussion also holds, {\em mutatis
mutandis,} for the general linear group, the orthogonal group, etc.}\
\ Two symplectic vector orbi-spaces are isomorphic if there exists a
linear symplectic isomorphism $\beta : V \to V'$ such that
$\beta(\Gamma)\beta^{-1} = \Gamma'$.

Let $\SP(V/\Gamma)$ denote the group of isomorphisms of $V/\Gamma$;
it is the group $N(\Gamma)/\Gamma$, where $N(\Gamma)$ is the normalizer
of $\Gamma$ in $\SP (V)$.  A {\em symplectic representation
} of a group $H$ on a symplectic vector
orbi-space $V/\Gamma$ is a group homomorphism $\rho : H \to \SP(V/\Gamma)$.
Two symplectic representations $\rho : H \to \SP(V/\Gamma)$
and $\rho' : H \to \SP(V'/{\Gamma'})$ are isomorphic if there
exists an isomorphism $\beta: V/\Gamma \to V'/{\Gamma'}$ such that
$\rho = \beta \rho' \beta^{-1}$.
In particular, two symplectic representations $\rho : H \to \SP(V/\Gamma)$
and $\rho' : H \to \SP(V/{\Gamma})$ are isomorphic if there
exists  $\beta \in \SP(V/\Gamma)$ such that $\rho = \beta \rho' \beta^{-1}$.
\begin{lemma}\label{lem_orbi-rep}
Let $\rho : H \to \SP (V/\Gamma)$ be a symplectic representation of a
group $H$ on a symplectic vector orbi-space $V/\Gamma$, and let
$N(\Gamma)$ denote the normalizer of $\Gamma$ in $\SP (V)$.  The
representation $\rho$ and the short
exact sequence $1\to \Gamma \to N(\Gamma) \to \SP (V/\Gamma) \to 1$
give rise to the {\em pull-back extension} $\pi: \hat{H} \to H$
and the {\em (symplectic) pull-back  representation}
$\hat{\rho}:\hat{H} \to N(\Gamma) \subset \SP (V)$
so that $\Gamma$ is naturally a subset of $\hat{H}$,
and the following diagram is exact and commutes.
\begin{equation}\label{eq_orbi-rep}
	\begin{diagram}
\node{1}\arrow{e} \node{\Gamma } \arrow{e}\arrow{s,=,-} \node{\hat{H}}
\arrow{s,r}{\hat{\rho}} \arrow{e,t}{\pi} \node{H} \arrow{s,r}{\rho} \arrow{e}
\node{1}\\
\node{1}\arrow{e} \node{\Gamma } \arrow{e} \node{N(\Gamma)} \arrow{e}
\node{\SP (V/\Gamma)} \arrow{e}\node{1}
	\end{diagram}
\end{equation}
If $\rho$ is faithful then $\hat{\rho}$ is also faithful.

Conversely, given a symplectic representation $ \hat{\rho}:\hat{H} \to
\SP (V)$ of a group $\hat{H}$ on a symplectic vector space $V$ and a
finite normal subgroup $\Gamma$ of $\hat{H}$, such that $\hat{\rho}$
is faithful on $\Gamma$, there exist a symplectic orbi-representation
$\rho : H \to \SP (V/\Gamma)$ of the quotient $H = \hat{H}/\Gamma$
making the diagram~\eqref{eq_orbi-rep} commute.
\end{lemma}
\begin{Proof}
Pull-backs exist in the category of groups.
\end{Proof}
\begin{lemma}\label{lem_lin_momentmap}
Let $\rho : H \to \SP (V/\Gamma)$ be a symplectic representation of a
group $H$ on a symplectic vector orbi-space $(V/\Gamma, \omega)$.
The action of $H$ on $V/\Gamma$ is Hamiltonian with
a moment map $\phi _{V/\Gamma}: V/\Gamma \to \fh^*$ given by the formula
\begin{equation}
	\label{eq_lin_moment}
 \langle \xi, \phi_{V /\Gamma}(v) \rangle  =
 \frac{1}{2} \omega(\xi \cdot v, v) \quad
 \text{ for all } \xi \in \hat{\fh} \text{ and } v\in V,
\end{equation}
where $\xi \cdot v$ is the value at $v$ of the vector field on $V/\Gamma$
induced by the infinitesimal action of $\xi \in \fh $.

The diagram
$$
\begin{diagram}
\node{V}\arrow{e}\arrow{s,l}{\hat{\phi}_V} \node{V/\Gamma }
\arrow{s,r}{\phi_{V/\Gamma}} \\
\node{\hat{\fh}^*}\node{\fh^* } \arrow{w,t}{\varpi ^*}
\end{diagram}
$$
commutes,  where $\hat{\phi}_V$ is the moment map for the action on $V$
of the pull-back extension $\pi :\hat{H}\to H$, and $\varpi
:\hat{\fh} \to \fh$ is the isomorphism of Lie algebras induced by the
homomorphism $\pi :\hat{H} \to H$.
\end{lemma}

\begin{proof}
It is easy to see that equation~\eqref{eq_lin_moment} defines a moment
map.  To show that the diagram commutes, it is enough to show that the
moment map $\hat{\phi} _V :V \to \hat{\fh}^* $ for the action of the
pull-back extension $\hat{H}$ on $V$ is $\Gamma$ invariant.  But
$\Gamma$ commutes with the identity component of $\hat{H}$.
\end{proof}

The following lemma is folklore.  A proof is given in Appendix~\ref{appendix}.
\begin{lemma}
\label{cor weight rep}
There is a bijective correspondence between isomorphism classes
of $2n$ dimensional symplectic representations
of a torus $H$ and unordered $n$ tuples of  elements
(possibly with repetition) of the weight lattice $\ell^* \subset \fh^*$ of $H$.

Let $(V,\omega)$ be a $2n$ dimensional symplectic vector space.
Let $\rho : H \to \SP (V, \omega )$ be a symplectic representation
with weights  $(\beta_1, \ldots, \beta_n)$.
There exists a decomposition $(V, \omega) =
\oplus _i (V_i, \omega _i)$ into invariant mutually perpendicular
2-dimensional symplectic subspaces and an invariant norm $|\cdot |$
compatible with the symplectic form $\omega = \oplus \omega _i$ so
that the representation of $H$ on $(V_i, \omega _i)$ has weight $\beta_i$
and
the moment map $\phi _\rho :V \to \fh^*$ is given by
\begin{equation} \label{eq moment map}
\phi _\rho (v_1, \ldots, v_n) =  \sum |v_i|^2 \beta_i \quad
\text{ for all } v= (v_1, \ldots, v_n) \in \oplus _i V_i.
\end{equation}
\end{lemma}

\begin{corollary}\label{cor rational cone}
Let $\rho : H \to \SP (V/\Gamma)$ be a symplectic representation of a
{\em torus} $H$ on a symplectic vector orbi-space $(V/\Gamma,
\omega)$.  The image of the corresponding moment map
$\phi_{V/\Gamma}(V/\Gamma)$ is a polyhedral cone in $\fh^*$ which is
rational with respect to the lattice ${\ell}\subset \fh$, the kernel
of the exponential map of $H$, a.k.a.\ the lattice of circle subgroups
of $H$.
\end{corollary}
\begin{proof}
This follows from Lemma~\ref{lem_lin_momentmap} and equation~\eqref{eq
moment map} in Lemma~\ref{cor weight rep}.
\end{proof}

These linear symplectic actions on symplectic vector orbi-spaces are
not only relatively easy to understand; they lie at the heart of every
Hamiltonian action.  Given a compact Lie group $G$ which acts on a
symplectic orbifold $(M,\omega)$ in a Hamiltonian fashion, we can
define the {\em symplectic slice} at a point $x \in M$.  The
$2$-form $\omega$ induces a non-degenerate antisymmetric bilinear form
on $\tT_xM$.  Let $\tT (G\cdot x)^\omega $ be the symplectic
perpendicular to the tangent space of $G \cdot x$ with respect to this
form.  The quotient
$$
V = \tT (G\cdot x)^{\omega }/( \tT (G\cdot x)\cap \tT (G\cdot x)^{\omega }),
$$
is naturally a symplectic vector space, The structure group of $x$,
$\Gamma$, acts symplectically on $\tT_xM$ and acts trivially on $\tT
(G \cdot x)$; therefore, $\Gamma$ acts symplectically on $V$.  The
{\em symplectic slice } at $x$ is the symplectic vector orbi-space
$V/\Gamma.$ The linear action of $G_x$ on $\tT_xM/\Gamma$ is
symplectic and preserves $\tT(G \cdot x)$.  Therefore, it induces a
symplectic representation of $G_x$ on $V/\Gamma$, the {\em
(symplectic) slice representation}.

As in the case of manifolds, the {\em differential} slice at $x$ is
isomorphic, as a $G_x$ representation, to the product
$$
	\fg_x^\circ  \times V/\Gamma,
$$
where $\fg_x^\circ$ denotes the annihilator of $\fg_x$ in $\fg^*$.
Thus, by Proposition~\ref{slice_prop}, a neighborhood of the orbit
$G \cdot x$ in $(M, \omega)$ is equivariantly diffeomorphic to a
neighborhood of the zero section in the associated orbi-bundle
$$
	Y = G\times_{G_x} \left(\fg_x^\circ \times V/\Gamma \right).
$$

\begin{Lemma}
	\label{locsymp}
Let $G$ be a compact Lie group.  Let $G\cdot x$ be an {\em isotropic}
orbit in a Hamiltonian $G$ orbifold $(M, \omega)$ and let $G_x \to \SP
(V/\Gamma)$ be the symplectic slice representation at $x$.  For every
$G_x$ equivariant projection $A: \fg \to \fg_x$, there is a $G$
invariant symplectic form $\omega _Y$ on the orbifold $Y= G\times_{G_x}
\left(\fg_x^\circ \times V/\Gamma \right) $ such that
\begin{enumerate}
\item
a neighborhood of  $G\cdot x$ in $M$
is equivariantly symplectomorphic to
a neighborhood of the zero section in $Y$,
and
\item
the action of $G$ on $(Y, \omega _Y)$ is Hamiltonian with a moment map
$\Phi_Y : Y \to \fg^*$ is given by
$$ \Phi _Y ([g, \eta, [v]])= Ad^\dagger (g) (\eta + A^* \phi
_{V/\Gamma} ([v])),
$$
\end{enumerate}
where $Ad^\dagger$ is the coadjoint action, $A^*: \fg^*_x \to \fg^* $
is dual to $A$, $\fg_x^\circ$ is the annihilator of $\fg_x$ in
$\fg^*$, and
$\phi_{V/\Gamma} : V/\Gamma \to \fg_x^*$ is the moment map for the
slice representation of $G_x$.
\end{Lemma}

\begin{Proof}
The construction is standard in the case of manifolds (cf.\
\cite{GS}); we adapt it for orbifolds.  Let $\hat{G}_x$ be the
pull-back extension of the isotropy group $G_x$ (cf.\
Lemma~\ref{lem_orbi-rep}).  The group $G_x$ acts on $G$ by $g_x \cdot
g = g g_x^{-1}$; this lifts to a symplectic action on the cotangent
bundle $T^*G$.  The corresponding diagonal action of $\hat{G}_x$ on
$T^* G \times V$ is Hamiltonian.  The projection $A: \fg \to \fg_x$
defines a left $G$-invariant connection 1-form on the principal $G_x$
bundle $G \to G/G_x$, and thereby identifies $Y$ with the reduced
space at zero $(T^* G \times V)_0$, thus giving $Y$ a symplectic structure.
The $G$ moment map on $T^* G \times V$ descends to a moment map for
the action of $G$ on $Y$, giving the formula in $(2)$.  The proof that
the neighborhoods are equivariantly symplectomorphic reduces to a form
of the equivariant relative Darboux theorem; it is identical to the
proof in the case of manifolds (see Remark~\ref{orbiproof}).
\end{Proof}

\begin{Remark}\label{uniqueness remark}{\em
Symplectic slice representations classify neighborhoods of orbits in
the following sense.  Let a compact Lie group $G$ act on symplectic
orbifolds $(M,\omega)$ and $(M',\omega')$ in a Hamiltonian fashion
with moment maps $\phi$ and $\phi'$ respectively.  Let $G \cdot x
\subset M$ and $G\cdot x' \subset M'$ be isotropic orbits.

Clearly, if there exist neighborhoods of $U$ of $G\cdot x$
and $U'$ of $G\cdot x'$ and a $G$ equivariant symplectomorphism $\psi:
U \to U'$ such that $\psi(x) = x'$, then the stabilizer of $x$ and
$x'$ is the same group $H$, and the slice representations at $x$ and
$x'$ are isomorphic.

Conversely, if the stabilizer of $x$ and $x'$ is the same group $H$,
and the symplectic slice representations at $x$ and $x'$ are
isomorphic, then it follows from Lemma~\ref{locsymp} that there exist
neighborhoods of $U$ of $G\cdot x$ and $U'$ of $G\cdot x'$ and a $G$
equivariant symplectomorphism $\psi: U \to U'$ such that $\psi(x) =
x'$ and $\phi' \circ \psi = \phi + const$.  }
\end{Remark}

\noindent

\begin{Remark}\label{abelian remark}\label{lem.loc_convex}
{\em Suppose again that a group $G$ acts in a Hamiltonian fashion on a
symplectic orbifold $(M,\omega)$ with moment map $\phi :M \to \fg^*$.
Suppose further that the group $G$ is a {\sl torus }.  Then the
coadjoint action of $G$ is trivial and every orbit is isotropic.  It
follows from Lemma~\ref{locsymp} that given an orbit $G\cdot x \subset
M$ there exist an invariant neighborhood $U$ of the orbit in $M$, a
neighborhood ${\cal W}$ of the zero section in the associated bundle
$G\times _{G_x} (\fg_x^\circ \times V/\Gamma)$ (where $G_x$ is the
isotropy group of $x$, $\fg_x^\circ$ is the annihilator of its Lie
algebra in $\fg^*$ and $V/\Gamma$ is the symplectic slice at $x$) and
an equivariant diffeomorphism $\psi :{\cal W} \to U$ sending the zero
section to $G\cdot x$ such that
\begin{equation}
(\phi \circ \psi) ([g, \eta, [v]]) = \alpha + \eta + A^* (\phi
_{V/\Gamma} ([v])),
\end{equation}
where $\alpha = \phi (x)$ and $A^* : \fg_x^* \hookrightarrow \fg^*$ is an
inclusion with $\fg^* = \fg^\circ_x \oplus A^* (\fg_x^*)$.
Consequently
\begin{equation}\label{eq nbd image}
	\phi (U) = (\alpha + W_1) \times (C \cap W_2),
\end{equation}
where $W_1 \subset \fg_x^\circ$ is a neighborhood of $0$ and $W_2
\subset A^* (\fg_x^*)$ is a neighborhood of the vertex of a rational
polyhedral cone $C$ ( $C = \phi _{V/\Gamma} (V/\Gamma)$ is a rational
polyhedral cone by Corollary~\ref{cor rational cone}).  }
\end{Remark}

The following result is a consequence of Lemma~\ref{locsymp} above.

\begin{corollary}
If a subgroup $H \subset G$ is connected, then $M^H$, the set of
points which are fixed by $H$, is a symplectic suborbifold.
\end{corollary}

\begin{Lemma}
\label{reduction lemma}
Let a compact group $G$ act in a Hamiltonian fashion on a symplectic
orbifold $(M,\omega)$  with moment map $\phi: M \to \ft^*$.  For a
regular value $\alpha \in \fg^*$ of $\phi$ which is fixed by the
coadjoint action, the {\em reduced space} of $M$ at $\alpha$,
$M_\alpha = \phi^{-1}(\alpha)/G$, is a symplectic orbifold.
\end{Lemma}

\begin{proof}
Consider $x \in \phi^{-1}(\alpha)$.  Because $\alpha$ is regular,
$G_x$, the stabilizer of $x$, is finite.  Since $\alpha$ is fixed by
the adjoint action, $G \cdot x$ is an isotropic orbit.  Let $G_x \to
\SP (V/\Gamma)$ be the symplectic slice representation at $x$.  By
Lemma~\ref{locsymp}, there is a $G$ invariant symplectic form on the
orbifold $Y= G\times_{G_x} \left(\fg^* \times V/\Gamma \right) $ such
that a neighborhood $U$ of $G\cdot x$ in $M$ is equivariantly
symplectomorphic to a neighborhood of the zero section in $Y$, and the
moment map $\Phi_Y : Y \to \fg^*$ is given by $\Phi _Y ([g, \eta,
[v]])= \alpha + Ad^\dagger (g) (\eta),$ where $Ad^\dagger$ is the
coadjoint action.

It is easy to see that $\phi^{-1}(\alpha) \cap U$ is isomorphic to
$G\times_{G_x} V/\Gamma$.  Therefore, a neighborhood of $[x]$ in
$M_\alpha$ is isomorphic to $V/\hat{G}_x$, where $\hat{G}_x$ is the
extension of $G_x$ by $\Gamma$.
\end{proof}

\begin{remark}{\em

\begin{enumerate}
\item The assumptions that the group $G$ is compact and that $\alpha$
is fixed by the coadjoint action are not necessary.  It is enough to
assume that the action of $G$ on the underlying topological space
$|M|$ is proper and that the coadjoint orbit through $\alpha $ is
locally closed.

\item If additionally we drop the assumption $\alpha$ is a regular
value of the moment map, then the quotient $\phi ^{-1} (G\cdot \alpha )/G$
is a symplectic stratified space in the sense of \cite{SL}.
\end{enumerate}
A proof of these two assertions would take us too far afield, so we
only note that the argument in \cite{BL} carries over to the case of
orbifolds.
}\end{remark}

\section{Morse Theory}\label{section.Morse}

In this section, we extend Morse theory to orbifolds.\footnote{ Since
orbifolds are stratified spaces, for Morse functions with isolated
fixed points this is a special case of Morse theory on
stratified spaces \cite{MTSS}.  } We need Morse theory for the
following result, which we will prove in the first part of this section.

\begin{Lemma}
\label{onedim}
Let $M$ be a connected compact orbifold,
and let $f: M \to \R$ be a Bott-Morse function
such that no critical suborbifold has index $1$ or $\dim(M)-1$.
The orbifold $M_{(a,b)} = f^{-1}(a,b)$ is connected for all $a, b \in \R$.
\end{Lemma}

We will use this result in the next section to prove that the fibers
of a torus moment map are connected, and that the image of a compact
symplectic orbifold under a torus moment map is a convex polytope.  In
the second part of the section we discuss  the notion of Morse
polynomials for functions on orbifolds and prove the Morse inequalities.

Most of the basic definitions needed for Morse theory on orbifolds are
strictly analogous to their manifold counterparts.  Let $f: M \to \R$
be a smooth function on an orbifold $M$.  A critical point $x$ of $f$
is {\em non-degenerate} if the Hessian $H(f)_x$ of $f$ is
non-degenerate.  More generally, a critical suborbifold $F \subset M$
is non-degenerate if for every point $x \in F$, the null space of the
Hessian $H(f)_x$ is precisely the tangent space to $F$.  In this case,
the Hessian restricts to a non-degenerate quadratic form $H$ on the
normal bundle $E$ of $F$ in $M$.

A smooth function $f: M \to \R$ is {\em Bott-Morse} if the set of
critical points is the disjoint union of non-degenerate suborbifolds,
and is {\em Morse} if each of these suborbifolds is a (isolated)
point.

\begin{Lemma}
Let $M$ be a compact orbifold.  Let $f: M \to \R$ be a smooth
function.  Choose $a < b \in \R $ such that $[a,b]$ contains no
critical values.  The orbifolds $M^-_a := f^{-1}(-\infty,a)$ and
$M^-_b := f^{-1}(-\infty,b)$ are diffeomorphic and so are the orbifolds
$f^{-1}(a)$ and $f^{-1}(b)$.
\end{Lemma}
\begin{proof}
The usual proof still applies, i.e., the diffeomorphism is given by
flows of the (renormalized) gradient of $f$ with respect to a
Riemannian metric.
\end{proof}

For critical points, the situation is only slightly more complicated.

\begin{Lemma} {\em ({\bf Morse Lemma})}
Let $F$ be a non-degenerate critical suborbifold of a smooth function
$f$ on an orbifold $M$.  Let $H$ be the restriction of the Hessian to
normal $E$ of $F$.  There exists a homeomorphism $\psi$ from a
neighborhood of the zero section of $E$ to a neighborhood of $F$ such
that
$$H = f \circ \psi. $$
\end{Lemma}

\begin{proof}
The proof the Morse lemma in \cite{MW} can be generalized to
vector bundles by carrying out the construction fiber by fiber.
Moreover, the construction is functorial once a metric has been chosen on $M$.
(See Remark \ref{orbiproof}.)

Alternatively, one can apply the same reasoning to Palais' proof of
the Morse lemma (cf. \cite{Lang}).  In this case, the map  $\psi$
will be a diffeomorphism.
\end{proof}

Let $F$ be an non-degenerate critical suborbifold of a smooth
function $f$ on a compact orbifold $M$.
Let $H$ be the restriction of the Hessian to the
normal bundle $E$ of $F$.
The bundle $E$ splits as a direct sum of vector
orbi-bundles $E^-$ and $E^+$ corresponding to the negative and
positive spectrum of $H$.
(By Remark~\ref{orbiproof}, this splitting exists
because such a splitting exists and is natural in the manifold case,
once a metric is chosen.)
The {\em index of $F$} is the rank of $E^-$.

We need only the following
homological consequence of the above lemmas.

\begin{lemma}
\label{morsebott}
Let $M$ be a compact orbifold, and let $f: M \to \R$ be a Bott-Morse
function.  Choose an interval $[a,b] \subset \R$ which contains a
unique critical value $c$.  Let $F$ be the critical suborbifold such
that $f(F) = c$.  Let $E^-$ be the negative normal bundle, and let
$D=D_F$ denote a disc bundle of $E^-$ and let $S= S_F$ denote the
corresponding sphere bundle.  The relative cohomology
$H_*(M^-_{b},M^-_{a}) = H_*(D_F,S_F)$.  Moreover, the boundary map
from $H_q(M^-_{b},M^-_{a})$ to $H_{q-1}(M^-_{b})$ in the long exact
sequence of relative homology is the composition of the isomorphism of
$H_*(M^-_{b},M^-_{a})$ and $H_*(D_F,S_F)$, the boundary map from
$H_q(D_F,S_F)$ to $H_{q-1}(S_F)$, and the map on homology induced by
the ``inclusion'' map $j$ from $S_F$ to $M^-_{a}$.  That is, the
following diagram commutes:
$$
\begin{diagram}
\node{H_q(M^-_{b},M^-_{a})}\arrow{s,l}{\partial}
\arrow{e,t}{\sim} \node{H_q(D_F,S_F)}\arrow{s,r}{\partial}\\
\node{H_{q-1}(M^-_{a})}\node{H_{q-1}(S_F)}\arrow{w,t}{j_*}
\end{diagram}
$$
\end{lemma}

\begin{proof}
Again, the manifold proof (see for example \cite{Chang}) can be
adapted.  The orbifold $M^-_{b}$ has the homotopy type of the space
obtained by attaching the disk orbi-bundle $D_F$ to $M_{a }$ by a map
from $S_F$ to $f^{-1}(a)$.  The result then follows from excision.
\end{proof}

We now prove Lemma~\ref{onedim} with a sequence of lemmas.

\begin{Lemma}\label{lemmanotone}
Let $F$ be an orbifold, let $\pi: E \to F$ be a $\lambda$ dimensional real
vector orbi-bundle,  and let $D$ and $S$ be the corresponding disk and
sphere orbi-bundles with respect to some metric.
If $\lambda = 0$, then $H_0(S) = 0$, so
$H_0(D,S)  \neq 0$.
Moreover, if $\lambda > 1$, then $H_1(D,S) = 0$.
In either case, the boundary map
$H_1(D,S) \to H_0(S)$ is trivial.
\end{Lemma}

\begin{proof}
If $\lambda = 0$ then $S$ is empty, so the result is trivial.

Suppose that  $\lambda > 1$;
we wish to show that $H_1(D,S) = 0$.
By the long exact sequence in relative homology it is enough to show
that the maps $H_0 (S) \to H_0 (D)$ and $H_1 (S) \to H_1
(D)$ induced by inclusion are injective and surjective, respectively.
But this follows
from two facts: the fibers of $\pi: S \to F$ are path connected and
any path in the base $F$ can be lifted to a path in the sphere bundle
$S$.
\end{proof}

\begin{Lemma}
\label{Morsesurj}
Let $M$ be a connected compact orbifold, and let $f: M \to \R$ be a
Bott-Morse function with no critical suborbifold of index $1$.  Then
\begin{enumerate}
\item $M^-_a:= \{m \in M: f(m) <a\}$ is connected for all $a \in \R$, and
\item  if $M^-_a \neq \emptyset$, then $H_1(M^-_a) \to H_1(M)$ is
a surjection.
\end{enumerate}
\end{Lemma}

\begin{Proof}
Let $F \subset M$ be a critical suborbifold of $f$ index $\lambda$.
Let $D_F$ and $S_F$ be the disk and sphere bundles of the negative
orbi-bundle over $F$.  Let $a = f(F)$, and let $\epsilon > 0$
be small.  Assume, for simplicity, that no other critical
suborbifold maps to $a$.

By Lemma~\ref{lemmanotone}, the map $H_1(D_F,S_F) \to H_0(S_F)$ is
trivial.  Therefore, by Lemma~\ref{morsebott}, the map
$H_1(M^-_{a+\epsilon},M^-_{a-\epsilon})) \to H_0(M^-_{a-\epsilon})$ is
also trivial.  Thus by the long exact sequence in relative homology,
the following sequence is exact:
$$
0 \to H_0(M^-_{a-\epsilon}) \to H_0(M^-_{a+\epsilon}) \to
H_0(M^-_{a+\epsilon},M^-_{a-\epsilon}) \to 0
$$
Therefore,
$\dim H_0(M^-_{a + \epsilon}) \geq \dim H_0(M^-_{a - \epsilon})$.  Since
$M$ is connected, this proves that $M_a^-$ is connected for all $a \in \R$.

If $\lambda = 0$, then $\dim H_0(M^-_{a + \epsilon})$ is strictly
greater than $\dim H_0(M^-_{a - \epsilon})$.  Therefore, since $M$ is
connected, the minimum is the unique critical value of index $0$.  For
any other critical value $a$, $H_1(M^-_{a+\epsilon},M^-_{a-\epsilon})
= 0$, by Lemma~\ref{lemmanotone} and Lemma~\ref{morsebott}.  Thus, by
the long exact sequence in relative homology, the map
$H_1(M^-_{a-\epsilon}) \to H_1(M^-_{a+\epsilon})$ is a surjection.
\end{Proof}

\begin{proof}{\bf of lemma~\ref{onedim}}
We may assume that $a$ and $b$ are regular values and that $M_{(a,b)}
:= \{ m \in M \mid a < f(m) < b \} \}$ is not empty.
By Lemma~\ref{Morsesurj},
$H_1(M^-_a) \oplus H_1(M^+_b) \to H_1(M)$ is a surjection,
where $M^+_b :=  f^{-1}(b,\infty)$.
Therefore, by Mayer-Vietoris, the following sequence is exact:
$$ 0 \to H_0(M_{(a,b)}) \to H_0(M^-_a) \oplus H_0(M^+_b) \to H_0(M)
\to 0.$$ By Lemma~\ref{Morsesurj}, $M^-_a$ and $M^+_b$ are
connected. Therefore $M_{(a,b)}$ is connected.
\end{proof}

\subsection{Morse polynomials}

We conclude the section with a few words about Morse polynomials for
orbifolds.  These observations are not used in the rest of the paper.

Let $M$ be a compact orbifold.
The {\em Poincar\'{e} polynomial} $\cp(x)$ is defined by
$$
\cp (x) =\sum_{i=0}^\infty \dim H^i(M) x^i.
$$
Let $f: M \to \R$ be a Bott-Morse function.  For each critical set
$F$, let $E_F^-$ be the vector orbi-bundle corresponding to the negative
spectrum of the Hessian of $f$ along $F$.  Let $D_F$ denote a disc
bundle of $E_F^-$ with respect to any metric, and let $S_F$ denote the
corresponding sphere bundle.  The {\em Morse polynomial} $\cm(x)$ is
defined by
$$
\cm(x)=\sum_F \sum_{i=0}^\infty \dim H^i(D_F,S_F) x^i,
$$
where the sum is taken over all critical orbifolds $F$.

\begin{theorem}
Let  $M$ be a compact orbifold, and let $f: M \to \R$ be a
Bott-Morse function.
If $\cp$ is the Poincare polynomial of $M$,
and $\cm$ is the Morse polynomial of $f$, then
$$\cm(x)- \cp(x) = (1+x) Q(x),$$
where $Q$ is a polynomial with
nonnegative coefficients.
\end{theorem}

\begin{proof}
Just as for manifolds, this follows immideatly from
Lemma~\ref{morsebott} and a spectral sequence argument.
\end{proof}

However, there is one very important distinction between Morse theory
on manifolds and on orbifolds.  For simplicity consider the case of
isolated critical points.  For a manifold, $\dim H^i(D,S)$ is easy to
compute; over any field, it is one if $\dim D = i$, and zero
otherwise.  Over finite fields, $\dim H^i(D/\Gamma,S/\Gamma)$ may be
much more complicated.  However, it is easy to see that
$H_i(D/\Gamma,S/\Gamma;\R) = 0$ if $i \neq \dim D$; whereas if $i =
\dim D$, then $H_i(D/\Gamma,S/\Gamma;\R) = \R$ if $\Gamma$ preserves
the orientation of $D$, and is trivial otherwise.  Therefore, the
$i$'th coefficient of the Morse polynomial is the number of points $x$
of index $i$ such that the the orbifold structure group of $x$
preserves the orientation of the negative eigenspace of the Hessian.

\begin{corollary}
The number of critical points with index $i$ is greater than or equal
to the dimension of $H_i(M)$.
\end{corollary}

\begin{example}
{\em
Let $M$ be torus, stood on end (see Figure~1 below).  Let $f: M \to
\R$ be the height function.  Let $\Gamma = \Z/(2\Z)$ act on $M$ by
rotating it $180$ degrees around the vertical axis.
Although there are four critical point,  $\cm(x) = 1 + 2x + x^2$.}
\end{example}
\setlength{\unitlength}{0.008in}
\begin{picture}(441,400)(10,-50)
\put(81,185){\ellipse{80}{156}}
\put(80,184){\ellipse{160}{240}}
\path(440,44)(439,329)
\path(441.028,321.007)(439.000,329.000)(437.028,320.993)
\path(182,183)(216,183)
\path(208.000,181.000)(216.000,183.000)(208.000,185.000)
\path(342,185)(386,185)
\path(378.000,183.000)(386.000,185.000)(378.000,187.000)
\path(322,63)	(317.538,63.540)
	(313.694,64.045)
	(310.404,64.527)
	(307.603,64.998)
	(303.210,65.951)
	(300.000,67.000)

\path(300,67)	(295.689,69.213)
	(293.183,70.748)
	(290.604,72.447)
	(285.716,75.957)
	(282.000,79.000)

\path(282,79)	(279.609,81.365)
	(276.868,84.374)
	(273.917,87.832)
	(270.899,91.542)
	(267.956,95.310)
	(265.230,98.939)
	(262.864,102.235)
	(261.000,105.000)

\path(261,105)	(258.511,109.277)
	(257.081,111.960)
	(255.625,114.812)
	(254.215,117.688)
	(252.926,120.442)
	(251.000,125.000)

\path(251,125)	(249.666,129.365)
	(248.963,132.083)
	(248.267,134.956)
	(247.599,137.835)
	(246.984,140.569)
	(246.000,145.000)

\path(246,145)	(245.010,149.381)
	(244.408,152.088)
	(243.791,154.947)
	(243.200,157.812)
	(242.677,160.539)
	(242.000,165.000)

\path(242,165)	(241.710,168.978)
	(241.516,173.834)
	(241.453,176.500)
	(241.412,179.275)
	(241.391,182.124)
	(241.390,185.009)
	(241.409,187.893)
	(241.445,190.741)
	(241.499,193.515)
	(241.569,196.179)
	(241.757,201.030)
	(242.000,205.000)

\path(242,205)	(242.481,209.221)
	(242.850,211.812)
	(243.273,214.539)
	(243.724,217.260)
	(244.179,219.837)
	(245.000,224.000)

\path(245,224)	(246.144,228.680)
	(246.906,231.551)
	(247.741,234.570)
	(248.608,237.579)
	(249.468,240.425)
	(251.000,245.000)

\path(251,245)	(252.753,249.297)
	(253.899,251.913)
	(255.148,254.653)
	(256.441,257.373)
	(257.719,259.933)
	(260.000,264.000)

\path(260,264)	(263.219,268.468)
	(265.314,271.116)
	(267.578,273.857)
	(269.894,276.553)
	(272.147,279.067)
	(276.000,283.000)

\path(276,283)	(280.956,287.218)
	(284.088,289.681)
	(287.446,292.196)
	(290.873,294.624)
	(294.213,296.826)
	(297.307,298.664)
	(300.000,300.000)

\path(300,300)	(303.061,301.052)
	(307.253,302.007)
	(309.927,302.477)
	(313.069,302.958)
	(316.739,303.462)
	(321.000,304.000)

\path(321,263)	(316.651,261.327)
	(312.917,259.833)
	(309.738,258.489)
	(307.051,257.267)
	(302.909,255.067)
	(300.000,253.000)

\path(300,253)	(296.152,248.198)
	(294.346,245.330)
	(293.000,243.000)

\path(293,243)	(290.975,238.912)
	(289.825,236.355)
	(288.660,233.641)
	(287.538,230.910)
	(286.515,228.300)
	(285.000,224.000)

\path(285,224)	(283.891,219.626)
	(283.312,216.910)
	(282.745,214.041)
	(282.210,211.166)
	(281.727,208.435)
	(281.000,204.000)

\path(281,204)	(280.412,199.831)
	(280.081,197.254)
	(279.758,194.535)
	(279.464,191.811)
	(279.225,189.221)
	(279.000,185.000)

\path(279,185)	(279.215,180.557)
	(279.454,177.832)
	(279.749,174.964)
	(280.075,172.101)
	(280.408,169.389)
	(281.000,165.000)

\path(281,165)	(281.688,160.344)
	(282.138,157.477)
	(282.644,154.458)
	(283.194,151.445)
	(283.777,148.593)
	(285.000,144.000)

\path(285,144)	(286.902,139.622)
	(288.217,137.012)
	(289.662,134.296)
	(291.148,131.606)
	(292.591,129.073)
	(295.000,125.000)

\path(295,125)	(296.375,122.700)
	(298.158,119.822)
	(302.000,115.000)

\path(302,115)	(304.652,113.231)
	(308.406,111.433)
	(310.835,110.465)
	(313.707,109.419)
	(317.077,108.272)
	(321.000,107.000)

\path(321,303)	(319.271,298.589)
	(318.090,295.306)
	(317.000,291.000)

\path(317,291)	(316.893,287.638)
	(317.140,283.449)
	(317.567,279.285)
	(318.000,276.000)

\path(318,276)	(318.908,271.452)
	(319.771,267.866)
	(321.000,263.000)

\path(321,106)	(319.271,101.589)
	(318.090,98.306)
	(317.000,94.000)

\path(317,94)	(316.893,90.638)
	(317.140,86.449)
	(317.567,82.285)
	(318.000,79.000)

\path(318,79)	(318.908,74.451)
	(319.771,70.866)
	(321.000,66.000)

\put(193,196){\makebox(0,0)[lb]{$f$}}
\put(357,196){\makebox(0,0)[lb]{$f$}}
\put(51,30){\makebox(0,0)[lb]{$\text{torus}$}}
\put(279,36){\makebox(0,0)[lb]{$\text{torus}/\Z_2$}}
\put(450,312){\makebox(0,0)[lb]{$\R$}}
\put(200,0){\makebox(0,0)[lb]{{\rm Figure 1}}}
\end{picture}

\section{Connectedness and Convexity}

In this section, we prove an analogue of the Atiyah connectedness
theorem and the Atiyah-Guillemin-Sternberg convexity theorem for
Hamiltonian torus actions on symplectic orbifolds \cite{GS} \cite{A}.
That is, we prove that the fibers of the moment map are connected, and
that the image of the moment map is a rational convex polytope Our
proofs are similar to Atiyah's proofs (op. cit.).  The fibers are
connected because the components of the moment map are Bott-Morse
functions with even indices, and convexity is an consequence of
connectedness.  The precise statements of these theorems follow.

\begin{Theorem}
\label{connected}
Let a torus $T$ act on a compact connected symplectic orbifold
$(M,\omega)$, with a moment map $\phi:M \to \ft^*$.  For every $a \in
\ft^*$, the fiber $\phi^{-1}(a)$ is connected.
\end{Theorem}

\begin{Theorem}
\label{convex}
Let a torus $T$ act on a compact connected symplectic orbifold
$(M,\omega)$, with a moment map $\phi:M \to \ft^*$.  The image of the
moment map $\phi(M) \subset \ft^*$ is a rational convex polytope.  In
particular it is the convex hull of the image of the points in $M$
fixed by $T$,
$$
	\phi (M) = \text{convex hull } (\phi (M^T)).
$$
\end{Theorem}

To prove  the theorems above, we will use the following lemma.
\begin{lemma}
\label{lemma Bott-Morse}
Let a compact Lie group $G$ act on a compact symplectic orbifold
$(M, \omega)$ with moment map $\phi : M \to \fg^*$.
For every $\xi \in \fg$ the $\xi
^{\text{th}}$ component of the moment map $\phi ^\xi := \xi \circ
\phi$ is Bott-Morse and the index of every critical orbifold is even.
\end{lemma}

\begin{proof}
This is a generalization of Theorem~5.3 in \cite{GS} and of
Lemma~2.2 in \cite{A} to the case of orbifolds.  The proof is the
same, except one must use the orbifold version of the equivariant
Darboux theorem (cf. Lemma~\ref{locsymp} which specializes to the
equivariant Darboux theorem when the orbit is a point).
\end{proof}

\begin{Remark} {\em
Moreover, if a component of the moment map $\phi^\xi$ has isolated
critical points, the orbifold isotropy groups of the critical points
preserves the symplectic form, and hence the orientation, on the
negative eigenspace of the Hessian.  Therefore, these  maps are
perfect Morse functions, i.e., the $i^{\text{th}}$ coefficient of the
Morse polynomial equals the $i^{\text{th}}$ coefficient of the
Poincar\'{e} polynomial, $\cm_i = \dim H_i (M)$.}
\end{Remark}

\begin{Proof} {\bf \!\!of Theorem \ref{connected}}
We will prove that the preimage of any ball is connected by induction
on the dimension of $T$.  Because the moment map $\phi$ is continuous
and proper, this implies that the fibers of $\phi$ are connected.

Consider first the case of $\dim T =1$.
By Lemma~\ref{lemma Bott-Morse}, moment maps for circle actions are
Bott-Morse functions of even index.
By lemma~\ref{onedim}, the preimage of any ball is connected.

Suppose that $T$ is a $k$ dimensional torus, $k>1$ and let $B$ be a
closed ball in $\ft^*$.  Let $\ell \subset \ft$ denote the lattice of
circle subgroups of $T$.  For every $0\not =\xi \in \ell$ the map
$\phi ^\xi \equiv \xi \circ \phi$ is a moment map for the action of
the circle $S_\xi := \{ \exp t\xi : t\in \R\}$.  Let $\R_\xi $ denote
the set of regular values of $\phi ^\xi$.  For every $a\in \R_\xi$ the
reduced space $M_{a, \xi} := (\xi \circ \phi)^{-1} (a)/S_\xi$ is a
symplectic orbifold (Lemma~{reduction lemma}).  The $k-1$ dimensional
torus $H:= T/S_\xi$ acts on $M_{a, \xi}$ in a Hamiltonian fashion.

The affine hyperplane $\{\eta \in \ft^*: \xi (\eta ) = a\}$ is
naturally isomorphic to the dual of the Lie algebra of $\fh$.  This
isomorphism identifies the restriction of $\phi$ to $(\phi ^\xi)^{-1}
(a)$ with the pull-back of an $H$-moment map $\phi ^H$ by the orbit
map $\pi :(\phi ^\xi)^{-1} (a) \to M_{a, \xi}$.  By inductive
assumption the preimages of balls under $\phi^H$ are connected.
Therefore, $\phi ^{-1} (B \cap \{\eta : \xi (\eta ) = a\}) = \pi ^{-1}
( (\phi ^H)^{-1} (B))$ is connected.

Now the set
$$
	U = \bigcup _{\xi \in \ell} \bigcup _{a\in \R_\xi } B \cap \{\eta : \xi
(\eta ) = a\}
$$
is connected and dense in the ball $B$, and its preimage
$\phi^{-1}(U)$ is connected.  Therefore the closure
$\overline{\phi^{-1} (U)}$ in $M$ is connected.  By
Lemma~\ref{locsymp}), $\overline{\phi^{-1} (U)} =
\phi^{-1}(\overline{U})$. Hence $\phi
^{-1}(B)=\phi^{-1}(\overline{U})$ is connected.
\end{Proof}

\begin{Proof} {\bf \!\!of Theorem~\ref{convex} }
Without loss of generality the action of $T$ is
effective, and hence, by Corollary~\ref{cor.loc_free}, free on a dense
subset.   Consequently the interior of  $\phi (M)$ is nonempty.

To prove that $\phi (M) $ is convex it suffices to show that the
intersection ${\cal L}\cap \phi (M)$ is connected for any rational
affine line ${\cal L}\subset \ft^*$, i.e., any line of the form ${\cal
L} = \R \upsilon +a$ for some $a$ in $\ft^* $ and $\upsilon$ in the
weight lattice $\ell ^*$ of $T$.  Define $\fh := \ker \upsilon$.  Let
$i^*$ be the dual of the inclusion $i:\fh \to \ft$, and let $\alpha =
i^*(a)$.  The map $\phi ^H := i^* \circ \phi : M \to \fh^*$ is the
moment map for a subtorus $H= \exp (\fh)$ of $T$.  The fibers $(\phi
^H)^{-1} (\alpha )$ are connected by Theorem~\ref{connected}.  On the
other hand,
$$
 (\phi ^H)^{-1} (\alpha) = \phi ^{-1} ((i^*)^{-1}(\alpha)) = \phi ^{-1} (\phi
(M) \cap (a + \R \upsilon)).
$$

Next, we show that $\phi (M)$ is the convex hull of $\phi (M^T))$.  By
Minkowski's theorem, since the set $\phi (M)$ is compact and convex,
it is the convex hull of its extreme points.  Recall that a point
$\alpha $ in the convex set $A$ is {\em extreme} for $A$ if it {\em
cannot} be written of the form $\alpha = \lambda \beta
+(1-\lambda)\gamma$ for any $\beta, \gamma \in A$ and $\lambda \in
(0,1)$. It follows from Lemma~\ref{locsymp} and Remark~\ref{abelian
remark} that for any point $m$ in a Hamiltonian $T$ orbifold $M$ the
image $\phi (M)$ contains an open ball in the affine plane $\phi (m) +
\ft_m^\circ$, where $\ft_m^\circ$ is the annihilator of the isotropy
Lie algebra of $m$ in $\ft^*$.  Therefore the preimage of extreme
points of $\phi (M)$ consists entirely of fixed points.

To show that $\phi (M)$ is a convex polytope it suffices to show that
$\phi (M^T)$ is finite.  This follows from that facts that $M^T$ is
closed, $M$ is compact, and $\phi$ is locally constant on $M^T$.
Therefore, we can write $\phi(M)$ as
$$
\phi(M) = \cap_{i = 1}^N \{ \alpha \in \ft^* \mid \left< \alpha, \xi_i
\right> \leq \eta_i \mbox{ for all } i\},
$$
for some $\eta_i \in \R$ and $\xi_i \in \ft$, where $N$ is the number
of facets.

To prove that $\phi (M)$ is rational we need to show that for each
$\xi_i$, the subgroup $H = H_i \subset T$, which is the closure of
$\{\exp t\xi_i: t\in \R\}$, is a circle.  The global maximum of the
function $\phi ^{\xi_i}$ is $\eta_i$.  Therefore the points in
$(\phi^{\xi_i}) ^{-1}(\eta_i)$ are fixed by $H$.  On the other hand,
for a generic $m$ in $(\phi^{\xi_i}) ^{-1}(\eta_i)$, the image
$d\phi(T_mM) \subset \ft^*$, is codimension one, so the stabilizer of
$m$ is one dimensional.
\end{Proof}

\part{Toric varieties}

\section{Local models for symplectic toric orbifolds}

In this section, we construct local models for symplectic toric
orbifolds. We begin with linear actions.

\begin{lemma}
	\label{linear}
Let $\rho: H \to \SP (V/\Gamma)$ be a faithful representation of a
compact abelian Lie group $H$ on a symplectic vector orbi-space such that
$\dim V/\Gamma = 2 \dim H$. The pull-back extension $\pi: \hat{H}\to H$
(cf.\ Lemma~\ref{lem_orbi-rep}) is a torus.

A fortiori, $\Gamma$ is abelian, and $H$ is a torus.
\end{lemma}

\begin{proof}
Since the orbi-representation $\rho : H \to \SP(V/\Gamma)$ is
faithful, the pull-back representation $\hat{\rho} : \hat{H} \to
\SP(V)$ is also faithful. Let $\tilde{H}$ denote the identity
component of $\hat{H}$.  The subgroup $\hat{\rho} (\tilde{H})$ is an $n$
dimension torus in $\SP (V)$, where $n = \frac{1}{2} \dim V $.  Every
maximal compact subgroup of $Sp(V)$ is isomorphic to $U(n)$. Every
maximal torus of $U(n)$ is $n$ dimensional.  Therefore, since the only
elements of $U(n)$ which commute with a maximal torus
are the elements of that torus, we will be done once we show that every
$a \in \hat{H}$ commutes with $\tilde{H}$.

Given $a \in \hat{H}$, define a continuous homomorphism $f_a: \hat{H}
\to \hat{H}$ by $f(h) = a h a^{-1}$.  Notice that $f_a (e) = e$, where
$e$ is the identity element of $\hat{H}$.  Moreover,
since $H$ is abelian,
$\pi(f_a (h)) = \pi(a) \pi(h) \pi(a)^{-1} = \pi(h)$.
Therefore,
for all $h \in \tilde{H}$, $f_a (h) = h$, i.e., $a$ commutes with
$\tilde{H}$.
\end{proof}

Let $\fh$ be a vector space with a lattice $\ell \subset \fh$.
A vector $g \in \ell$ is {\em primitive} if
there is no positive integer $n>1$ such that  $\frac{1}{n}g \in \ell$.
A {\em rational cone} in $\fh^*$ is a set $C$ such that
$$
	C = \{ \eta \in \fh^* \mid \langle f_i, \eta \rangle \geq 0\},
$$
where $f_i \in \ell$ for all $i$.  We may assume the normals $f_i$ are
primitive.  A cone is {\em strictly convex} if it contains no
nontrivial subspace.  A cone is {\em simplicial} if it spans $\fh^*$
and has $\dim \fh$ facets.

\begin{lemma}\label{lemma weights-cones}
Let $\rho: H \to \SP(V/\Gamma)$ be a faithful symplectic
representation of a torus $H$ on a symplectic vector orbi-space such that
$\dim V = 2 \dim H$, and let $\phi_\rho: V/\Gamma \to \fh^*$ be the
corresponding moment map.  The image $C = \phi_{\rho}(V/\Gamma)$ is a
rational strictly convex simplicial cone (relative to the
lattice of circle subgroups $\ell \subset \fh$).
For every facet $F$ of this cone there
exists a positive integer $m_F$ such that $\Z/m_F\Z$ is the orbifold
structure group of every point in $\phi_\rho^{-1}(\intF)$,
the preimage of the interior of $F$.

Conversely, given a a strictly convex simplicial rational cone $C$ in
$\fh^*$ and a set of positive integer $\{m_F\}$ indexed by the facets
of $C$, there exists a unique faithful symplectic representation $\rho
 :H \to\SP( V/\Gamma)$ such that $\phi_\rho (V/\Gamma ) = C$ and such that
 for every facet $F$ of $C$ the orbifold structure group of every point in
the preimage $\phi_\rho^{-1}(\intF)$ is $\Z/m_F\Z$.
\end{lemma}

\begin{proof}
We first show how an orbi-representation gives rise to a labeled
cone.

Let $\pi: \hat{H} \to H$ be the the pull-back extension, let
$\hat{\rho}$ be the pullback representation of $\hat{H}$ on $V$, and
let $\phi _{\hat{\rho}} :V \to \hat{\fh}$ be associated moment map.
By Lemma~\ref{cor weight rep}, the vector space $V$ can be split into
the direct sum of symplectic, invariant two dimensional subspaces, $V
= \oplus V_i$.  Let $\hat{\rho}_i :\hat{H}\to V_i$ denote the $i$'th
subrepresentation, and let $f_i \in \hat{\ell}^*$ denote the
corresponding weight.

Since the representation $\hat{\rho}$ is faithful, the vectors
$\{f_i\}$ form a basis of the weight lattice.  Let $\{e_i\}$ denote
the dual basis of the lattice $\hat{\ell}$.  By equation
\eqref{eq moment map} in Lemma~\ref{cor weight rep} the image $\phi
_{\hat\rho}(V)$ is the strictly convex simplicial rational cone.
$$
\phi _{\hat\rho}(V) =
	 \{ \eta \in \hat\fh^* \mid \langle e_i, \eta \rangle \geq 0\}.
$$

By Lemma~\ref{lem_lin_momentmap}, the diagram
$$
\begin{diagram}
\node{V}\arrow{e}\arrow{s,l}{\hat{\phi}_V} \node{V/\Gamma }
\arrow{s,r}{\phi_{V/\Gamma}} \\
\node{\hat{\fh}^*}\node{\fh^* } \arrow{w,t}{\varpi ^*}
\end{diagram}
$$
commutes, where $\varpi :\hat{\fh} \to \fh$ is the isomorphism of Lie
algebras induced by $\pi: \hat{H} \to H$, and $\varpi^*$ is its
transpose.  Therefore, image $\phi _\rho (V/\Gamma)$ is the rational
simplicial cone
$$
	C = \{ \eta \in \fh^* \mid \langle \varpi (e_i), \eta \rangle \geq 0\},
$$
the open facets $\intF_j$ of $C$ are given by
$$
\intF_j = \{ \eta \in \fh^* \mid \langle \varpi (e_i), \eta \rangle>
0\, \text{ for } i\not =j,\, \langle \varpi (e_j), \eta \rangle = 0\}
$$
and their preimages are
$$
\phi _\rho ^{-1} (\intF_j) = \{ [(v_1, \ldots, v_n)] \in (\oplus
V_i)/\Gamma \mid v_j=0,\, v_i \not =0 \text{ for } i\not =j\}.
$$
Consequently for any $[v] \in \phi _\rho ^{-1} (\intF_j) $ the
stabilizer of $v \in V$ in $\hat{H}$ is $\cap_{i \neq j}
\ker(\hat{\rho_i})$, the circle subgroup with Lie algebra $\R e_j$.
Therefore, the stabilizer of $v$ in $\Gamma$ is $\Z/m_j \Z$, where
$m_j$ is defined by fact that the primitive element of the
intersection of the ray $\bbR_{\geq 0} \varpi (e_j)$ with the lattice
$\ell$ is $\frac{1}{m_j} \varpi (e_j)$.

Conversely, suppose we are given a rational strictly convex simplicial
cone $C$ in $\fh^*$, and a collection of positive integers $\{m_F\}$
indexed by the set ${\cal F} (C)$ of facets of $C$.  For every facet
$F\in {\cal F} (C)$ there exists a primitive vector $g_F \in \ell$
such that
$$
C = \{ \eta \in \fh^* \mid \langle g_F \eta \rangle \geq 0, \, \text{
for all }\, F \in{\cal F} (C) \}.
$$
Because $C$ is simplicial, the set
$\{e_F = m_F g_F\}_{F\in {\cal F} (C)}$ is a basis of a sublattice
$\hat{\ell}$ of $\ell$.
Let $\{f_F = e_F^*\}$ denote the dual basis of $\hat{\ell}^*$.
Let $\hat{H}$ be the torus $ \fh/\hat{\ell}$.
By Lemma~\ref{cor weight rep} there exists a unique symplectic
reperesentation $\hat{\rho} : \hat{H} \to \SP(V)$ with weights $\{f_F
\}$. Clearly $\dim V = 2 \dim \hat{H}$.  Since $\{f_F \}$ is a basis,
the representation $\hat{\rho}$ is faithful.

Let $\Gamma = \ell/\hat{\ell}$.  There is a short exact sequence $1
\to \Gamma \to \hat{H} \stackrel{\pi}{\to} H \to 1$.  Thus
$\hat{\rho}$ defines an orbi-representation $\rho : H \to \SP (V/\Gamma)$.
We leave it for the reader to show that the image of the corresponding
moment map $\phi _\rho$ is the cone $C$ and that the orbifold
structure groups of points in the preimages of open facets are cyclic
groups of correct orders.
\end{proof}

\begin{lemma}\label{local1} \label{local}
Let $(M,\omega,T, \phi )$ be a compact symplectic toric orbifold.
\begin{enumerate}
\item
The moment map is an orbit map, i.e., for any $a\in \phi (M)$ the
fiber $\phi ^{-1} (a)$ is a single $T$ orbit.  Moreover, for any $a\in
\phi (M)$ there exists a neighborhood $W_a$ in $\ft^*$ such that $\phi
^{-1} (W_a)$ is a tubular neighborhood of $\phi ^{-1} (a)$.

\item \label{item connected}
The isotropy group of every point $x\in M$ is connected.
\end{enumerate}
\end{lemma}

\begin{proof} Let $x$ be a point in the fiber $\phi ^{-1} (a)$.
By Remark~\ref{abelian remark} there exist an invariant neighborhood
$U$ of the orbit $T\cdot x$ in $M$, a neighborhood ${\cal W}$ of the
zero section in the associated bundle $T\times _{T_x} (\ft_x^\circ
\times V/\Gamma)$ (where $T_x$ is the isotropy group of $x$,
$\ft_x^\circ$ is the annihilator of its Lie algebra in $\ft^*$ and
$V/\Gamma$ is the symplectic slice at $x$) and an equivariant
diffeomorphism $\psi :{\cal W} \to U$ sending the zero section to
$T\cdot x$ such that
\begin{equation}
(\phi \circ \psi) ([g, \eta, [v]]) = \alpha + \eta + A^* (\phi
_{V/\Gamma} ([v])),
\end{equation}
where $\alpha = \phi (x)$ and $A^* : \ft_x^* \hookrightarrow \ft^*$ is
an inclusion such that $\ft^* = \ft^\circ_x \oplus A^* (\ft_x^*)$.

Since $\dim M = 2 \dim T$ and the action of $T$ on $M$ is effective,
$\dim V/\Gamma = 2 \dim H$ and the
action of $H$ on $V/\Gamma$ is effective.
These two facts, together with Lemma~\ref{linear}, imply that $T_x$
is connected.

The same two facts imply that ${\phi_{V/\Gamma}}^{-1}(0)$ is a single
point (cf.\ proof of Lemma~\ref{lemma weights-cones}).  Consequently
$\phi^{-1}(a) \cap U = T \cdot x$.  Since $\phi^{-1}(a)$ is connected
by Theorem~\ref{connected} and since the orbit $T \cdot x$ is closed,
it follows that $T \cdot x = \phi^{-1}(a)$.  Therefore, $\phi$ is an
orbit map.

Since $\phi$ is proper, given a neighborhood $U$ of the fiber $
\phi^{-1}(a)$ there exists a neighborhood  $W$ of $a$ in $\ft^*$ such
that $\phi^{-1}(W) \subset U$.
\end{proof}

Recall that by Theorem~\ref{convex}, given a compact symplectic toric
orbifold $(M,\omega,T, \phi )$, the image $\phi (M)$ is a polytope.
Combining Lemma~\ref{local1} with Lemma~\ref{locsymp} we get the
following description of neighborhoods of fibers of moment maps on
compact symplectic toric orbifolds.

\begin{corollary}\label{toric orbifold model}
Let $(M,\omega,T, \phi )$ be a compact symplectic toric orbifold, let
$F$ be a face of the polytope $\phi(M)$
and let $a$ be a point in its interior.

Let $H$ be the torus with Lie algebra $\fh$,
where the affine plane spanned by $F$
is of the form $\fh^\circ +a$,
and $\fh^\circ$ is the annihilator of $\fh$.
For any point $x\in \phi^{-1} (a)$, the isotropy group of $x$ is $H$.

Let $K$ be a complimentary subtorus so that $T\simeq K \times H$.  For a
sufficiently small neighborhood $W_a$ of $a$ in $\ft^*$ there exists
an open equivariant symplectic embedding $\psi: \phi ^{-1}
(W_a)\hookrightarrow T^* K \times V/\Gamma$, where $V/\Gamma$ is the
symplectic slice at $x$.  In ``coordinates''defined by the embedding
$\psi$, the moment map $\phi \circ \psi^{-1} : \psi (\phi ^{-1}
(W_a))\to \fk^*\times \fh^* \simeq \ft^* $ is given by
\begin{equation}\label{moment equation 1}
	(k,\eta, [v]) \mapsto (a+ \eta, \phi_{V/\Gamma} ([v]),
\end{equation}
where $\phi_{V/\Gamma} :V/\Gamma \to \fh^*$ is the moment map for the
slice representation.  Moreover, the action of $H$ on $V/\Gamma$ is
faithful and consequently $C := \phi_{V/\Gamma}(V/\Gamma)$ is a simplicial
cone.  Finally,
\begin{equation}\label{moment equation 2}
   W_a \cap \phi (M) = W_a \cap \left( (a + \fh^\circ) \times C \right).
\end{equation}
\end{corollary}

\begin{corollary}\label{simple polytope}
Let $(M,\omega,T, \phi )$ be a compact symplectic toric orbifold.
The set $\phi (M)$ is a {\em simple} rational polytope.
\end{corollary}
\begin{proof}
Immediate from Corollary~\ref{toric orbifold model}
\end{proof}

\begin{lemma}\label{lemma weights}
Let $(M,\omega,T, \phi )$ be a compact symplectic toric orbifold and
let $F$ be a facet of $\phi (M)$.  For any $x$ in $\phi ^{-1}(\intF)$,
the preimage of the interior of $F$, the orbifold structure group
$\Gamma_x$ is cyclic, and the order of the group $\Gamma_x$ depends
only on the facet $F$.
\end{lemma}
\begin{proof}
Since the facet $F$ has codimension 1, it follows from
Corollary~\ref{toric orbifold model} that the isotropy group $H$ of
$x$ is a circle and that $V/\Gamma$, the symplectic slice at $x$, is 2
dimensional.  Consequently the orbifold structure group $ \Gamma _x$
is cyclic.  Moreover, it follows from \eqref{moment equation 1} that
for a sufficiently small neighborhood $W_a$ of $a$ the preimage $\phi
^{-1}(W_a \cap F)$ is isomorphic to $K\times\{0\}\times \{0\}$ in
$T^*K \times V/\Gamma_x$, where $T^*K$ is the cotangent bundle of the
subtorus complimentary to $H$.  Therefore the function that assigns to
a point in the interior $\intF$ of $F$ the order of the orbifold
structure group of a point in its preimage is locally constant on
$\intF$. Since $\intF$ is connected, the result follows.
\end{proof}

Putting together Lemma~\ref{lemma weights} with Corollary~\ref{simple
polytope} we get the following result.

\begin{corollary}
	\label{labeled polytope}
Let $(M,\omega,T, \phi )$ be a compact symplectic toric orbifold.  Let
${\cal F}(\phi (M))$ denote the set of facets of $\phi (M)$ and for a
facet $F$ of $\phi (M)$ let $m_F$ be the order of the structure group
$\Gamma _x$ for some $x\in \phi ^{-1} (\intF)$.  Then $(\phi (M),
\{m_F\}_{ F\in {\cal F}(\phi (M))})$ is a labeled polytope in the
sense of Definition~\ref{def two}.
\end{corollary}
The labeled polytope $(\phi (M), \{m_F\}_{ F\in {\cal F}(\phi (M))})$
associated by Corollary~\ref{labeled polytope} to a compact
symplectic toric orbifold $(M,\omega,T, \phi )$ determines the
orbifold locally in the following sense.

\begin{lemma}
\label{uniqueness of local models}
Let $(M,\omega,T, \phi )$ be a compact symplectic toric orbifold and
let $(\phi (M), \{m_F\}_{ F\in {\cal F}(\phi (M))})$ be the associated
labeled polytope.  For a small enough neighborhood $W_a$ of a point
$a\in (\phi (M))$ the symplectic toric orbifold $(\phi ^{-1} (W_a),
\omega|_{(\phi ^{-1} (W_a)}, T, \phi)$ is uniquely determined by the
sets $W_a$ and $\{(F, m_F) \mid F\in {\cal F}(\phi (M)), \, a\in F\}$.
\end{lemma}
\begin{proof}
 Let $x$ be a point in the fiber $\phi ^{-1} (a)$. By
Lemmas~\ref{local1} and~\ref{locsymp} the orbifold $\phi ^{-1} (W_a)$
is uniquely determined by the symplectic slice representation $\rho :
H \to \SP(V/\Gamma )$ at $x$.  To determine the isotropy group $H$ of
$x$ we observe that the intersection $L = \cap \{F\in {\cal F}(\phi
(M)) \mid a\in F\}$ is the face of the polytope $\phi (M)$ containing
$a$ in its interior, and that by Corollary~\ref{toric orbifold model}
the face $L$ determines the isotropy group of $a$.

By Lemma~\ref{lemma weights-cones}, in order to determine the
representation itself, it is enough to determine the image of
corresponding moment map and the orders of the relevant orbifold
structure groups. The orders of the structure groups are the integers
$\{m_F \mid a\in F\}$.  Let $K$ be a subtorus complementary to $H$ in
$T$. The splitting $\ft = \fk \times \fh$ defines a projection $p$
from $\ft^*$ onto $\fh^*$ along the affine subspace $\fh^\circ +a
$. The projection maps the intersection of a small neighborhood $W_a$
of $a$ with $\phi (M)$ onto a neighborhood of a vertex of a simplicial
rational strictly convex cone $C \subset \fh^*$.  The pair $(C, \{m_F
\mid a\in F\})$ determines the slice representation $H \to \SP
(V/\Gamma)$.
\end{proof}

\section{From Local to Global}\label{local to global}

In this section, which is joint work with Chris Woodward, we show that two
compact symplectic toric orbifolds which have isomorphic associated labeled
polytopes
are themselves isomorphic.  First, we have already shown
that they are {\em locally}
isomorphic.  By extending proposition~2.4 in \cite{HS} to the symplectic
category, we show that these local isomorphisms can be glued together to
construct a global isomorphism.

We say that two compact symplectic toric orbifolds $(M,\omega,T,\phi)$
and $(M',\omega',T,\phi')$ with $\phi(M) = \phi'(M') = \Delta$ are
{\em locally isomorphic over} $\Delta$ if every point in $\Delta$ has
a neighborhood $U$ such that ${\phi'}^{-1}(U)$ and $\phi^{-1} (U)$ are
isomorphic as symplectic toric oribifolds.

\begin{definition}
\label{def.sheaf}
Let $(M,\omega,T,\phi)$ be a compact symplectic toric orbifold; let $\Delta =
\phi(M)$.  Define a sheaf $\cS$ over $\Delta$ as follows: for each open set
$U \subset \Delta$, $\cS(U)$ is the set of isomorphisms of $\phi^{-1}(U)$.
\end{definition}

\begin{lemma}
\label{lem.sheaf}
Let $(M,\omega,T,\phi)$ be a compact symplectic toric orbifold.
Define a sheaf $\cS$ over $\Delta = \phi(M)$ as in Definition
\ref{def.sheaf}.
The cohomology group $H^1(\Delta,\cS)$ classifies (up to an isomorphism)
compact symplectic toric orbifolds $(M',\omega',T,\phi')$ such that $\phi'(M')
=
\Delta$ and $M'$ is locally isomorphic to $M$ over $\Delta$.
\end{lemma}

\begin{proof}
Let $\U = \{ U_i \}_{i \in I}$ be a covering of $\Delta$ such that there
is an isomorphism $h_i : \phi^{-1}(U_i) \to {\phi'}^{-1}(U_i)$ for
each $i \in I$.  Define $f_{ij}: \phi^{-1}(U_i \cap U_j) \to
\phi^{-1}(U_i \cap U_j)$ by $f_{ij} =  h_i^{-1} \circ h_j$.
The set $\{f_{ij}\}$'s is a closed cochain in $C^1(\U,\cS)$.  Moreover,
the cohomology class of this cocycle is independent of the choices of
the isomorphisms $h_i$.

Conversely, given a cocycle $\{f_{ij}\} \in C^1(\U,\cS)$, we can construct
a compact symplectic toric orbifold with moment polytope $\Delta$ by taking the
disjoint union of the $\phi^{-1}(U_i)$'s and gluing $\phi^{-1}(U_i)$ and
$\phi^{-1}(U_j)$ together using the isomorphisms $f_{ij}$.
\end{proof}

\begin{proposition}
\label{loctoglob}
Let $(M,\omega,T,\phi)$ be a compact symplectic toric orbifold.
Define a sheaf $\cS$ over $\Delta = \phi(M)$ as in Definition
\ref{def.sheaf}.
The sheaf $\cS$ is abelian and
the cohomology group $H^i(\Delta,\cS)$ is  $0$ for all $i > 0$.
\end{proposition}

\begin{proof}
Let $\underline{\ell \times \R}$ denote the sheaf of
locally constant functions with values in $\ell \times \R$, where
$\ell \subset \ft$ is the lattice of circle subgroups.
Since $\Delta$ is contractable, $H^i(\Delta,\underline{\ell
\times \R}) = 0$ for all $i > 0$.

Define a sheaf $\cC^\infty$ over $\Delta$ as follows: for each open set $U
\subset \Delta$, $\cC^\infty(U)$ is the set of smooth $T$ invariant
functions on $\phi^{-1}(U)$.  We may think of elements of $\cC^\infty(U)$
as continuous functions on $U$ which pull back to smooth functions on
$\phi^{-1}(U)$.
\footnote{
One can show using the Lemma~\ref{local}
and a theorem of G.W. Schwarz \cite{sc:sm} that for every $T$
invariant smooth function $f$ on $\phi^{-1}(U)$ there exists a smooth
function $\varphi_f $ on $\ft^*$ with $f =\varphi_f \circ \phi $ on
$\phi^{-1}(U)$.} \,

Since $\cC^{\infty}$ is a fine sheaf
$H^i(\Delta, C^{\infty}) = 0$ for all $i > 0$.

Therefore,  to prove that $\cS$ is abelian and that $H^i(\cS,\Delta) = 0$
for $i > 0$, it suffices to construct the following sequence of sheaves,
and to show that it is exact:
$$
	0 \to  \underline{\ell \times \R} \stackrel{j}{\to}
\cC^\infty \stackrel{\Lambda}{\to} \cS \to 0.
$$

Define $j : \underline{\ell \times \R} \to \cC^\infty$ as follows:
given $(\xi,c) \in \ell \times \R$ and a point $x \in M$,
let $j(\xi,c)(x) = c +  \left<\xi,\phi(x)\right>.$

Next, we construct the map $\Lambda: C^\infty \to \cS$.  Let $U \subset
\Delta$ be an open set and let $f : \phi^{-1}(U) \to \R$ be a smooth $T$
invariant function.  The flow of the Hamiltonian vector field $\Xi_f$ of $f$
on $\phi^{-1}(U)$ is $T$ equivariant and preserves the moment map $\phi$.
Define $\Lambda(f)$ to be the time one flow of $\Xi_f$.
We now show that the sequence of sheafs is exact.\\[-8pt]

The map $j$ is clearly injective.\\[-8pt]

Recall that for every vector $\xi \in \ft$ there exists a vector field
$\xi _M$ on $M$ induced by the action of $T$, that $\xi _M$ is the
Hamiltonian vector field of the function $\left<\xi,\phi(x)\right>$.
The time $t$ flow of $\xi _M$ is given by
\begin{equation}
	\label{flow}
	x\mapsto e^{t\xi} \cdot x ,
\end{equation}
where $\xi \mapsto e^\xi$ is the exponential map from $\ft$ to $T$.
Therefore $\im j \subset \ker \Lambda$.

To show that $\im j \supset \ker \Lambda$ we argue as follows.  Without
loss of generality we may assume that the subset $U$ of $\Delta$ is the
intersection of $\Delta $ with a ball in $\ft^*$. Let $U_0$ be the
intersection of this ball with the interior $\intD$ of the polytope
$\Delta$. Then both $U$ and $U_0$ are convex hence contractible.  By
Lemma~\ref{toric orbifold model}, $\phi ^{-1}(U_0)$ is open and dense in $\phi
^{-1}(U)$ and $\phi : \phi ^{-1}(U_0) \to U_0$ is a principal $T$ bundle.

Let $f: \phi ^{-1}(U) \to \R$ be a $T$ invariant function with $\Lambda (f)
= id$.  We want to show that $df = d \langle \xi, \phi\rangle$ on $\phi
^{-1}(U)$ for some $\xi \in \ell$.  It is enough to show that this equality
holds on $ \phi ^{-1}(U_0)$.

Since $f$ is $T$ invariant and since $\phi : \phi ^{-1}(U_0) \to U_0$ is a
principal $T$ bundle, there exists $h\in C^\infty (U_0)$ with $f = h\circ
\phi$.  Hence the Hamiltonian vector field $\Xi_f$ of $f$ at the points $x\in
\phi ^{-1}(U_0)$ is given by
$$
	\Xi_f (x) = (dh (\phi (x))_M (x)
$$
(since $dh (\phi (x)) \in T^* _{\phi (x)} \ft^* = \ft$, the expression $(dh
(\phi (x))_M $ makes sense).

Equation (\ref{flow}) implies that if $X: \phi ^{-1} (U) \to \ft$ is a $T$
invariant function and $Y$ is a vector field on $\phi ^{-1} (U)$ defined by
$Y(x) = (X(x))_M (x)$ then the time $t$ flow $\psi_t$ of $Y$ is given by
\begin{equation}
	\label{flow2}
	\psi_t : x\mapsto e^{tX(x)} \cdot x .
\end{equation}
Consequently the time one flow $\Lambda (f)$ of $\Xi_f$ is given by
$$
	\Lambda (f): x\mapsto e^{dh(\phi(x))} \cdot x,
		\qquad x\in  \phi ^{-1}(U_0).
$$
Since by assumption $\Lambda (f) (x) = x$ for all $x\in \phi ^{-1}(U_0)$,
we have $dh(u) \in \ell$ for all $u\in U_0$.  Since $U_0$ is connected and
since $\ell$ is discrete, the continuous function $dh: U_0 \to \ell$ is
constant.  Thus $df = d \langle \xi, \phi\rangle$ for some $\xi \in \ell$
and all $x\in \phi ^{-1}(U_0) $.\\[-3pt]

The final step is to show that $\Lambda$ is surjective. If the ball used in
the definition of the set $U$ is small enough then by Lemma~\ref{local} the
set $\phi^{-1}(U)$ is a tubular neighborhood of some $T$ orbit in
$\phi^{-1}(U)$.  Let $\psi$ be an isomorphism of $\phi^{-1}(U)$.  We must
show there exists a $T$ invariant function on $\phi^{-1}(U)$ whose time one
flow is the map $\psi$.

Since $\psi$ is an isomorphism, it is, a fortiori, a $T$-equivariant
diffeomorphism of $\phi^{-1}(U)$ which preserves orbits.  Therefore, by
Theorem~3.1 in \cite{HS}, there exists a smooth $T$ invariant map $\sigma :
\phi^{-1}(U) \to T$ such that $\psi(x) = \sigma (x)\cdot x$.  Since $U$ is
contractible and $\sigma $ is $T$ invariant, there exists a smooth map
$X:\phi^{-1}(U) \to \ft$ such that $e^{X (x)} = \sigma (x)$.  As before
define a vector field $Y$ on $\phi^{-1}(U)$ by $Y(x) = (X(x))_M (x)$.
By equation (\ref{flow2}) the time one flow of $Y$ is $x\mapsto e^{X(x)}
\cdot x = \sigma (x)\cdot x = \psi (x)$.  Thus it is enough to show
that $Y$ is a Hamiltonian vector field, i.e., that the contraction
$\iota (Y) \omega $ is exact.

Just as for a free action of a compact Lie group on a manifold, we
can, following Koszul \cite{Kosz}, define on the orbifold $M$ a
{\em complex} of basic forms.  Namely, a form $\alpha \in \Omega (M)$ is
{\em basic} if $\alpha $ is $T$ invariant and if for any vector $\xi
\in \ft$, we have $\iota (\xi _M) \alpha =0$.  Similarly we can define
basic forms on any open $T$ invariant subset of $M$, such as $\phi
^{-1} (U)$.  We observe that basic forms have two properties.

\begin{enumerate}

\item A basic form $\alpha$ is preserved by any $T$ equivariant map $\psi
:M \to M$ which induces the identity map on the orbit space $M/T$, that is,
$\psi ^* \alpha = \alpha$.  This is true because it is a closed condition,
which holds on the open dense smooth subset of the orbifold $M$ where the
action is free.

\item The integral of a basic $k$ form over a $k$ cycle which lies entirely
in a $T$ orbit is zero.  It follows that the cohomology of the complex
of basic forms on a tubular neighborhood of an orbit is trivial.  In
other words, a closed basic form is exact on a tubular neighborhood of
an orbit.
\end{enumerate}

We now argue that the contraction $\iota (Y) \omega $ is a closed
basic form on the tubular neighborhood $\phi ^{-1} (U)$.  Then, by
property $(2)$ above, there exists a basic zero form $f$ such that
$\iota (Y) \omega = df$.

Since $Y$ and $\omega$ are $T$ invariant, $\iota (Y) \omega $ is $T$
invariant.  Since the $T$ orbits are isotropic in $M$, and since $Y$ is
tangent to $T$ orbits, $\iota (\xi _M) \iota (Y) \omega = 0$ for any $\xi
\in \ft$.  The Lie derivative $L_Y \omega = d \iota (Y) \omega $ is also
basic, since basic forms are a subcomplex.  Consequently, since the time
$t$ flow $\psi _t$ of $Y$ induces the identity map on the orbit space, we
have by property $(1)$ above that $\psi_t ^* L_Y \omega = L_Y \omega$.
Since $\psi _1 = \psi $ and since $\psi $ is symplectic, we have
\begin{eqnarray*}
0= \psi _1 ^* \omega - \omega &=& \int _0^1 \frac{d}{dt} \psi _t^*
\omega\, dt\\
&= & \int _0^1 \psi _t^* L_Y \omega\, dt\\
&= & \int _0^1 L_Y \omega\, dt = L_Y \omega = d  \iota (Y) \omega .
\end{eqnarray*}
This proves there exists a $T$ invariant function $f$ whose time 1
flow is the isomorphism $\psi$. Hence $\Lambda$ is surjective.
\end{proof}

\begin{Theorem} \label{uniqueness}
Two compact symplectic toric orbifolds which have isomorphic associated labeled
polytopes are themselves isomorphic.
\end{Theorem}

\begin{proof}
Without loss of generality, we may assume that the labeled polytopes
associated to two orbifolds $M$ and $M'$ are equal.  By
Lemma~\ref{uniqueness of local models}, $M$ and $M'$ are locally
isomorphic.  By Proposition~\ref{loctoglob}, $H^1(\Delta,\cS) = 0.$ By
Lemma~\ref{lem.sheaf}, this implies that $M$ and $M'$ are isomorphic
as symplectic toric orbifolds.
\end{proof}

\begin{remark}{\em
Given any labeled polytope $\Delta$, one can construct local models for
the symplectic toric orbifold associated to $\Delta$.  Since we've shown
that $H^2(\cS,\Delta) = 0$, an argument in \cite{HS} similar to Lemma
\ref{lem.sheaf} shows that there exists a symplectic toric orbifold which
corresponds to the given labeled polytope.  In section~\ref{surj}, we give
a more explicit construction.  }
\end{remark}

\section{Existence}
\label{surj}

Given any labeled polytope, we construct a compact symplectic toric
orbifold such that its associated labeled polytope is the one which we
began with.  This construction is a slight variation of Delzant's
construction \cite{Del}.

\begin{theorem}\label{existence}
Let $T$ be a torus.  Let $\ft$ denote its Lie algebra, and let
$\ell \subset \ft$ denote the lattice of circle subgroups.
Given a simple rational polytope $\Delta \subset \ft$
and a positive integer $m_F$ attached to each facet $F$ of $\Delta$,
there exists a compact symplectic toric orbifold $(M, \omega, T, \phi)$  such
that
$\phi (M) = \Delta$ and the orbifold structure group of a point
in $M$ which maps to the interior of a facet $F$ is $\bbZ/ m_F \bbZ$.

Moreover, $(M,\omega,T,\phi)$ is a symplectic reduction of $\C^N$ by
an abelian subgroup of $\Su(N)$.
\end{theorem}

\begin{proof}
The polytope $\Delta$ can be written uniquely as
$$
\Delta = \cap_{i = 1}^N
\{ \beta  \in \ft^* \mid \left<  \beta , m_i y_i \right> \leq \eta_i \},
$$
where $N$ is the number of facets, the vector $y_i \in \ell$ is the primitive
normal to the $i$th facet, $m_i$ is the integer attached to the $i$th
facet, and $\eta = (\eta_1, \cdots, \eta _N) \in (\bbR ^N)^*$.

Define a linear projection $\varpi : \R^N \to \ft$ by $\varpi(e_i) =
m_i y_i$, where $\{e_i\}$ is the standard basis for $\R^N$.  This
defines a short exact sequence and its dual:
$$
0 \lra \fk  \stackrel{j}{\lra} \R^N \stackrel{\varpi}{\lra} \ft \lra 0
$$
and
$$
0 \lra \ft^* \stackrel{\varpi^*}{\lra} (\R^N )^* \stackrel{j^*}{\lra} \fk^*
\lra 0.
$$
Let $K$ denote the kernel of the map from
$\bbT ^N= \R^N/\Z^N $ to $T = \ft/\ell$
which is induced by $\varpi$. The kernel is given by
$$
  K= \{ [\alpha ]\in R^N/\Z^N \mid \sum_{i=1}^N \alpha_i m_i y_i \in \ell \}
$$
The Lie algebra of $K$ is $\fk$, the kernel of $\varpi$.

Consider $\C^N$ with the standard symplectic form $\sum \sqrt{-1} dz_j
\wedge d\bar{z}_j$.  The standard ${\bbT}^N$ action on $\C^N$ has
moment map $\phi_{{\bbT}^N}(z_1,\ldots,z_N) = \sum _{i=1}^N |z_i|^2
e_i^* = (|z_1|^2, \ldots, |z_N|^2)$, where $\{e^*_i\}$ is the basis
dual to $\{e_i\}$.  Clearly, $\left<\phi_{{\bbT}^N}(z),e_i\right> \geq
0$ for all $i$, and $\left<\phi_{{\bbT}^N}(z),e_i\right> = 0$ exactly
if $z_i = 0$.  Since $K$ is a subgroup of ${\bbT}^N$, $K$ acts on
$\C^N$ with moment map $\phi_K = j^* \circ \phi_{{\bbT}^N}$.

The stabilizer of a point $z \in \C^N$ in $\bbT^N$ is
$$
\bbT^N_z = \{ [\alpha] \in \bbR^N/\bbZ^N \mid e^{2\pi \sqrt{-1} \alpha_i} = 1
\mbox{ for all } i \mbox{ with } z_i \neq 0 \}.
$$
Since the stabilizer of $z \in \C^N$ in $K$ is  $\bbT^N_z \cap K$,
\begin{equation*}
\label{Kstab}
K_z =
\{ [\alpha] \in \bbR^N/\bbZ^N \mid e^{2\pi
\sqrt{-1} \alpha_i} = 1 \mbox{ for all } i \mbox{ with } z_i \neq 0 \mbox{
and } \sum_{i=1}^N \alpha_i m_i y_i \in \ell \}.
\end{equation*}

Define an affine embedding $\iota_{\eta} : \ft^* \to ({\R^N})^*$ by
$\iota _\eta(\beta ) = \varpi^*(\beta ) - \eta$.  Note that
\begin{eqnarray*}
\iota_\eta (\Delta)& = &\{ \xi \in (\R^N)^* \mid \xi \in \varpi ^* (\ft^*)
-\eta \mbox{ and } \xi_i \geq 0 \mbox{ for all } i \}\\
& = &\{ \xi \in (\R^N)^* \mid \xi \in (j^*)^{-1} ( j^*(-\eta)) \mbox{ and }
\xi_i \geq 0 \mbox{ for all } i \}
\end{eqnarray*}
Moreover,
$$
 \iota _\eta (\Delta \cap \{\beta \in \ft^* \mid \langle \beta , m_j y_j
\rangle = \eta _j \}) = \iota_\eta(\Delta) \cap \{ \xi \in (\bbR
^N)^* \mid \xi _j = 0\}.
$$

For every $z\in \bbC ^N$ such that $\phi _K (z) = j^*(-\eta )$,
$\phi_{\bbT^N}(z)$ is in $\iota_\eta(\ft^*).$
Since $\Delta $ is simple, for every point $\beta \in \Delta$ the set
$$
	\{ y_i \mid \langle \beta , m_i y_i \rangle = \eta _i \}
$$
is linearly independent.
Consequently the set
$$
	\{ y_i  \mid z_i=0\}  = \{ y_i  \mid \phi_{\bbT^N}(z)_i=0\}
$$
is linearly independent.
Hence the isotropy group $K_z$ is discrete.  Therefore, $j^*(-\eta)$
is a regular value of $\phi_K$, and the reduced space
$M = \phi _K ^{-1 }(j^*(-\eta )) /K$
is a symplectic toric orbifold.

Since the action of $\bbT^N$ on $\bbC^N$ commutes with the action of $K$,
it induces a Hamiltonian action of $\bbT^N$ on $M$. Moreover, the moment
map $\phi _{\bbT^N}$ descends to a moment map $\tilde{\phi}: M\to
(\bbR^N)^*$ and $\tilde{\phi} (M) = \iota _\eta (\Delta )$.  In fact the
action of $\bbT^N$ on $M$ descends to a Hamiltonian action of $T= \bbT^N /K$
and $\phi _T = (\iota _\eta |_\Delta )^{-1} \circ \tilde {\phi }$ is a
corresponding moment map.

We claim that the action of $T$ on $M$ is effective.  It suffices to show
that there exists a point $z\in \phi _K ^{-1} (j^* (-\eta ))$ so that its
isotropy group in $\bbT^N$ is trivial.  Such a point exists because the
isotropy group in $\bbT^N$ of any point $z\in \phi _{\bbT^N}^{-1} (\{ \xi
\in (\bbR^N)^* \mid \xi_i >0 \mbox{ for all } i \})$ is trivial, the embedding
$\iota _\eta $ maps the interior of the polytope $\Delta $ into the set $\{
\xi \in (\bbR^N)^* \mid \xi >0 \mbox{ for all } i \})$ and $\Delta $ has
non-empty interior.

It remains to show that the orbifold structure group of a point $[z]$ in
$M$ mapping to the interior of the facet cut out by the hyperplane $\{\beta
\in \ft^* \mid \langle \beta , m_j y_j\rangle = \eta _j\}$ is $\bbZ/ m_j\bbZ$.
But $[z]$
lies in the interior of this facet if and only if $\phi _{\bbT^N}(z) \in
\{\xi \in (\bbR^N)^* \mid \xi _j = 0 \mbox{ and } \xi_i
\not =0 \mbox{ for } i\not =
j\}$. For such a point $z$ the isotropy group $K_z$ is $\bbZ/ m_j \bbZ$.
\end{proof}

\section{Compatible complex structures} \label{Compatible complex structures}

In this section, we show that every compact symplectic toric orbifold
possesses a invariant complex structure which is compatible with the
symplectic form.  Moreover, suppose that two compact symplectic toric
orbifolds $(M,\omega,T,\phi)$ and $(M',\omega',T,\phi')$ are given
invariant complex structure which are compatible with their
symplectic forms.  They are equivariantly biholomorphic exactly if the
polytopes $\phi(M)$ and $\phi'(M')$ have the same fan.

\begin{Theorem}
Every compact symplectic toric orbifold possesses an invariant
complex structure which is compatible with its symplectic form, i.e.,
every such orbifold has an invariant  K\"ahler structure.
\end{Theorem}

\begin{proof}
Let $(M,\omega,T,\phi)$ be a compact symplectic toric orbifold.

By Theorem \ref{simple polytope}, $\Delta = \phi(M)$ is a simple
rational polytope, By Theorem~\ref{lemma weights}, for each facet $F$
of $\Delta$, there exists a positive integer $m_F$ such that
$\Z/m_F\Z$ is the orbifold structure group of every point in $M$ which
maps to the interior of $F$.

By Theorem \ref{existence}, there exists a compact symplectic toric
orbifold $(M',\omega',T,\phi')$ which is a symplectic reduction of
$\C^N$ by an abelian subgroup of $\Su(N)$, such that $\phi'(M) =
\Delta$ and such that for every facet $F$ of $\Delta$, $\Z/m_F\Z$ is
the orbifold structure group of every point in $M'$ which maps to the
interior of $F$.  Since $M'$ is the reduction of a K\"ahler manifold
by a group which preserves its K\"ahler structure, by Theorem~3.5 in
\cite{GS:gq} $M'$ possesses an equivariant K\"{a}hler structure which
is compatible with its symplectic form.

By Theorem \ref{uniqueness}, $M$ and $M'$ are equivariantly
symplectomorphic; therefore, $M$ inherits an an invariant K\"{a}hler
structure which is compatible with its symplectic form.
\end{proof}

Not only do all compact symplectic toric orbifolds admit compatible complex
structures, they are, in fact, algebraic varieties.

\begin{lemma}
\label{lemma variety}
Let $(M, \omega, T, \phi)$ be a compact symplectic toric orbifold and
let $J$ be a $T$-invariant complex structure on $M$ which is
compatible with the symplectic form $\omega$. Then $M$ has the
structure of an algebraic toric variety with the fan equal to the fan
defined by the polytope $\phi (M)$.
\end{lemma}

\begin{remark}{\em
If the class of the K\"ahler form $\omega$ in $H^2(M)$ is rational
then by the Kodaira-Baily embedding theorem \cite{Baily}, $M$ is a
projective algebraic variety.  The projective embedding provided by
the Kodaira-Baily theorem is equivariant with respect to the action of
the torus $T$ and so $M$ is a projective toric variety.  The class
$[\omega]$ is rational if and only if the edges of the polytope $\phi
(M)$ are rational vectors relative to the weight lattice of $T$.
Therefore, in order to prove that {\em all} toric orbifolds are
algebraic, we  argue differently.}
\end{remark}

An immediate consequence of Lemma~\ref{lemma variety} is the following
theorem.

\begin{theorem}
Let two compact symplectic toric orbifolds $(M,\omega,T,\phi)$ and
$(M',\omega',T,\phi')$ be given invariant complex structure which are
compatible with their symplectic forms.  The orbifolds are
equivariantly biholomorphic exactly if the fans defined by their
polytopes are equal.
\end{theorem}

\begin{remark}{\em
This theorem shows, in particular, that if a compact symplectic toric
orbifold admits two different compatible complex structures, they are
equivariantly biholomorphic.  In contrast, the K\"ahler structure is
not unique.  For instance, there are many $S^1$ invariant K\"ahler
structures on $S^2$.}
\end{remark}

\begin{remark}{\em
There are several reasons for the difference of the classification of
symplectic toric orbifolds and of algebraic toric varieties.

The first reason is that some toric varieties do not admit any
symplectic form.  These correspond to fans which do not come
from a polytope.  The second reason is that changing the cohomology
class of the symplectic form corresponds to changing the length of the
edges of the polytope. This information is lost in the algebraic
category.  Finally, the integers attached to the faces are lost.  It
is easy to see why.  Give $\C$ the standard K\"ahler structure, and
let $\Z/m\Z$ act on $\C$.  The orbifolds $\C$ and $\C/\Z/m\Z$ are not
diffeomorphic, but they are $\C^\times $ equivariantly biholomorphic.
  }\end{remark}

\begin{proof}{\bf of Lemma~\ref{lemma variety}}
Since the complex structure $J$ is $T$-invariant, the action of $T$ on
$M$ extends to the action of the complexification $T_\C$ of $T$.  The
action of $T_C$ on $M$ has a dense open orbit. Denote it by $T_C\cdot
m$.

The action of $T_\C$ on $M$ can be linearized near fixed points.  That
is, if $x\in M$ is fixed by $T$, there exist a $T_\C$-invariant
neighborhood $V$ of $0$ in $T_x M$, a $T_\C$ invariant neighborhood of
$x$ in $M$ and abiholomorphic map $f: Y \to U$ which is
$T_\C$-equivariant.  There are several ways to see that the
linearization exists.  For example, the linearization proof for group
actions on K\"ahler manifolds in \cite{LS} is natural and so, by
Remark~\ref{orbiproof}, translates into a proof in the orbifold case
(See also \cite{Sj}, where more general result --- a holomorphic slice
theorem --- is proved.  The reader is refered to \cite{Sj} especially
since \cite{LS} may never see the light of day).  Alternatively we can
appeal to the holomorphic slice theorem in \cite{HL} which holds for
K\"ahler {\em spaces}, hence in particular, for orbifolds.

Since the action of $T$ on $M$ and hence on $T_x M$ is faithful and
since $\dim T = \frac{1}{2}\dim M$, the neighborhood $V$ is all of
$T_xM$ and $T_xM$ is a toric variety.  The linearization map embeds
this variety into $M$.  The fan of this variety consists of a single
simplicial cone together with its faces.  Moreover, this is the cone dual
to the image of the moment map corresponding to the linear action of
$T$ on the tangent space $T_xM$.

If $y$ is another fixed point and $h: T_y M \to M$ another
linearization, then both images $f(T_xM)$ and $h(T_yM)$ must contain
the dense open orbit $T_\C \cdot m$ of $T_\C$.  It is not hard to see
that the transition map from $f^{-1} (T_\C \cdot m)$ to $h^{-1} (T_\C
\cdot m)$ is rational and, in fact, is the same as the map defined by
the intersection of the corresponding fans.  The lemma now follows.
\end{proof}

\appendix

\section{Symplectic weights}
\label{appendix}

In this section we show that weights for symplectic representations of
tori are well defined by proving the following lemma.

\begin{lemma}
	\label{cor weight rep2}
There is a bijective correspondence between isomorphism classes of
$2n$ dimensional symplectic representations of a torus $H$ and
unordered $n$-tuples of elements (possibly with repetition) of the
weight lattice $\ell^* \subset \fh^*$ of $H$.

Let $(V,\omega)$ be a $2n$ dimensional symplectic vector space.
Let $\rho : H \to \SP (V, \omega )$ be a symplectic representation
with weights  $(\beta_1, \ldots, \beta_n)$.
There exists a decomposition $(V, \omega) =
\oplus _i (V_i, \omega _i)$ into invariant mutually perpendicular
2-dimensional symplectic subspaces and an invariant norm $|\cdot |$
compatible with the symplectic form $\omega = \oplus \omega _i$ so
that the representation of $H$ on $(V_i, \omega _i)$ has weight $\beta_i$
and
the moment map $\phi _\rho :V \to \fh^*$ is given by
\begin{equation} \label{eq moment map2}
\phi _\rho (v_1, \ldots, v_n) =  \sum |v_i|^2 \beta_i \quad
\text{ for all } v= (v_1, \ldots, v_n) \in \oplus _i V_i.
\end{equation},
\end{lemma}

\begin{proof}
Since symplectic representations do not have  naturally defined complex
structures, we define a {\em character} of a torus $H$ to be a
homomorphism into an {\em oriented } circle.
A {\em weight} is the differential of a character.

A symplectic representation $\rho : H \to \SP (V, \omega)$ on a
symplectic 2-plane $(V, \omega)$ has a well defined character.
Assuming that $\rho$ is non-trivial, the image $\rho (H)$ is a compact
abelian subgroup of $\SP (V, \omega)$, and hence is a circle.  Circle
subgroups of $\SP (V, \omega)$ are naturally oriented by the
symplectic form $\omega$ via their orbits in $V$.  Additionally, any
two circle subgroups of $\SP (V, \omega)$ are conjugate by an element
of $\SP (V, \omega)$, and conjugation preserves the induced
orientations.  We conclude that a conjugacy class of symplectic
representation of a torus on a symplectic 2-plane has a well-defined
character (and hence a well-defined weight) which determines this
conjugacy class uniquely.

Therefore, we may assume that the plane is
$\bbC$ with symplectic form $\sqrt{-1}dz\wedge d \bar{z}$ and that the
action of a torus $H$ is given by $(a,z) \mapsto e^{i \beta (\log a)}z$
where $\beta$ is a weight of $H$.  In agreement with formula~(\ref{eq moment
map2}), the moment map for this action is given
by $z\mapsto |z|^2 \beta$ and the image of the plane is a ray $\bbR_{\geq
0} \beta$ through the weight $\beta$.

By contrast, for the underlying real representation the weight is
defined only up to a sign.  For example, conjugation by
$\left(
\begin{array}{cc} 0& 1\\ 1&0 \end{array}\right)$
sends the matrix
$\left( \begin{array}{rr} \cos \theta & \sin \theta\\ -\sin \theta&
\cos \theta \end{array}\right)$
to
$\left( \begin{array}{rr} \cos
\theta & -\sin \theta\\ \sin \theta& \cos \theta \end{array}\right)$.

Since a maximal compact subgroup of $\SP (\bbR^{2n})\simeq \SP (V,
\omega)$ is isomorphic to $U(n)$, there exists on $V$ an $H$ invariant
complex structure $J$ which is compatible with the symplectic form,
i.e., $\omega (J\cdot ,J\cdot ) = \omega (\cdot, \cdot)$ and $g(\cdot,
\cdot) = \omega (\cdot, ,J\cdot )$ is positive definite.  Then $V$
decomposes as a direct sum of mutually orthogonal invariant complex
line. $(V, \omega) = \oplus _{i=1}^n (V_i, \omega _i)$.  These $V_i$'s
are mutually symplectically orthogonal invariant $2$ planes.  The
moment map is given by formula (\ref{eq moment map2}), where
$(\beta_1, \ldots, \beta_n)$ are the corresponding weights.

Although this decomposition is not natural, we will show that if a
symplectic representation $(V, \omega)$ of a torus $H$ has two
different decompositions, $(V, \omega) = \oplus _{i=1}^n (V_i, \omega
_i) = \oplus _{i=1}^n (V'_i, \omega' _i)$ with weights $(\beta_1,
\ldots, \beta_n)$ and $(\beta'_1, \ldots, \beta'_n)$ then the two
$n$-tuples of weights are the same up to a permutation.  For a weight
$\alpha \in \ell^*$ of the torus $H$ the isotypical subspace $W_\alpha
= \oplus _{\beta_i = \pm \alpha}V_i$ is canonically defined for the
underlying real representation.  Therefore, if a weight $\beta_i$
occurs in the first decomposition of $(V, \omega)$, then the subspace
$W_{\beta_i} = \oplus _{\beta'_j = \pm \beta_i} V'_j$ is nonempty. The
image of $W_{\beta_i}$ under the moment map is the Minkowski sum of
the rays through the weights $\beta'_j$ such that $\beta'_j= \pm
\beta_i$.  On the other hand, neither the vector space $W_{\beta_i}$
nor its image under the moment map depend on the decomposition of $V$
into symplectic planes.  Since the image contain the ray through
$\beta_i$, there must exist $\beta'_j$ with $\beta_i = \beta'_j$.  We
can then split off the $2$ plane corresponding to $\beta_i$ and repeat
the above argument.
\end{proof}

\end{document}